\documentclass{amsart}
\usepackage[margin=1.2in]{geometry}
\usepackage[utf8]{inputenc}
\usepackage[T1]{fontenc}
\usepackage{graphicx}
\usepackage{amsfonts}
\usepackage{amssymb}
\usepackage{amsmath}
\usepackage{amsthm}
\usepackage{paralist}
\usepackage{nicefrac}
\usepackage{mathtools}
\usepackage{color}
\usepackage{verbatim}
\usepackage{abstract}
\usepackage{dsfont}

\usepackage[pdfborder={0 0 0}]{hyperref}

\newtheorem{satz}{Proposition}[section]
\newtheorem{theo}[satz]{Theorem}
\newtheorem{lemma}[satz]{Lemma}

\theoremstyle{definition}
\newtheorem{defi}[satz]{Definition}
\newtheorem{remark}[satz]{Remark}
\newtheorem{bsp}[satz]{Example}
\numberwithin{equation}{section}

\newcommand\ko{\hspace{.1em}}
\newcommand{\e}{\operatorname{e}}
\renewcommand{\Im}{\operatorname{Im}}

\newcommand{\tr}{\operatorname{tr}}
\newcommand{\R}{\mathbb{R}}
\newcommand{\C}{\mathbb{C}}
\newcommand{\Z}{\mathbb{Z}}

\newcommand{\N}{\mathbb{N}}
\renewcommand{\H}{\mathbb{H}}

\renewcommand{\theta}{\vartheta}
\renewcommand{\epsilon}{\varepsilon}

\newcommand{\I}{\mathcal{I}}
\newcommand{\ur}{U_r}
\newcommand{\us}{U_s}
\newcommand{\vr}{V_r}
\newcommand{\vs}{V_s}
\newcommand{\trans}{\mathsf{T}\hspace{-2pt}}

\newcommand{\eulerop}{\textbf{E}}
\newcommand{\eulerentry}{\textbf{E}_{ij}}
\newcommand{\laplace}{\boldsymbol{\Delta}}
\newcommand{\laplacian}{\operatorname{tr}\boldsymbol{\Delta}}

\newcommand{\partialmatrix}[1]{\frac{\partial}{\partial #1}}
\newcommand{\mymatrix}[2]{{#2}^{{#1}\times {#1}}}

\newcommand{\myintegral}{\int\limits_{\mathbb{R}^{m\times n}}}

\newcommand{\Aneg}{A^{-}}
\newcommand{\Apos}{A^{+}}
\newcommand{\Amaj}{M}

\newcommand{\mythetac}[2]{\vartheta_{{#1},{#2}}}

\newcommand\smod[1]{\ (\operatorname{mod} #1)}

\newcommand{\uvector}{\boldsymbol{u}}
\newcommand{\ucolumnvector}[1]{\boldsymbol{u_{#1}}}

\newcommand{\vcolumnvector}[1]{\boldsymbol{v_{#1}}}
\newcommand{\upos}{\boldsymbol{u^+}}
\newcommand{\uneg}{\boldsymbol{u^-}}
\newcommand{\ucolumnvectorpos}[1]{\boldsymbol{u^+_{#1}}}
\newcommand{\ucolumnvectorneg}[1]{\boldsymbol{u^-_{#1}}}
\newcommand{\Upos}{U^+}
\newcommand{\Uneg}{U^-}
\newcommand{\hvector}{\boldsymbol{h}}
\newcommand{\kvector}{\boldsymbol{k}}

\newcommand{\myhomofct}[1]{\mathcal{F}_{#1}^{m, n}}
\newcommand{\myhomopol}[1]{\mathcal{P}_{#1}^{m,n}}
\newcommand{\basis}[2]{\mathcal{B}_{#1}^{m,#2}}

\newenvironment{acknowledgments}{
   \abstract}{
  \endabstract
}
\makeatletter
\@namedef{subjclassname@2020}{\textup{2020} Mathematics Subject Classification}
\makeatother

\title{Siegel theta series for indefinite quadratic forms}
\author{Christina Roehrig}
\address{University of Cologne, Department of Mathematics and Computer Science, Division of Mathematics, Weyertal 86-90, 50931 Cologne, Germany}
\email{croehrig@math.uni-koeln.de}
\begin{document}
\maketitle
\begin{abstract}
The modular transformation behavior of theta series for indefinite quadratic forms is well understood in the case of elliptic modular forms due to Vignéras, who deduced that solving a differential equation of second order serves as a criterion for modularity.
In this paper, we will give a generalization of this result to Siegel theta series.
\end{abstract}
\section{Introduction}
In the course of his work on the Minkowski-Hasse principle for quadratic forms over the rationals, Siegel introduced a natural generalization of elliptic modular forms of higher genus $n$ \cite{Sie35}. Among those functions, nowadays called Siegel modular forms, Siegel theta series play a similarly important role as do theta series do in the context of elliptic modular forms.
In a recently published article by  Dittmann, Salvati Manni, and Scheithauer \cite{DMS19}, a basis of the space of Siegel cusp forms of degree 6 and weight 14 is given by harmonic Siegel theta series. By considering one of these basis elements, the authors deduce that the Kodaira dimension of the Siegel modular variety $\mathcal{A}_6=\operatorname{Sp}_{12}(\Z)\setminus \H_6$ is non-negative.

In order to give more examples of Siegel theta series and make use of the connection to various topics -- such as algebraic geometry and number theory -- it is desirable to give a general framework for the description of holomorphic and non-holomorphic Siegel theta series analogous to what is already known for elliptic theta series owing to the work of Vignéras \cite{Vig77}.
If theta series are built from functions that satisfy a certain second-order differential equation, the modularity of these series  immediately follows. For the (generalized) error functions, which are employed in the recent discussions of theta series for indefinite quadratic forms, this criterion is used to derive the modular transformation behavior of the emerging theta series. Namely, these are the results by Zwegers \cite{Zwe02} for quadratic forms of signature $(m-1,1)$, by Alexandrov, Banerjee, Manschot, and Pioline \cite{ABMP18} for signature $(m-2,2)$ and for arbitrary signature by Nazaroglu \cite{Naz18} and Westerholt-Raum \cite{W17}, which are brought together by Kudla \cite{Kud18} and Funke and Kudla \cite{FK19}. 
Even before that, Kudla and Millson \cite{KM86,KM87} considered a certain class of Schwartz functions to define modular forms in terms of theta functions and obtain holomorphic modular forms valued in the set of cohomology classes.

In most of these examples the criterion given by Vignéras plays an important role in order to deduce modularity, so the question arises whether a similar result holds for more general types of theta series. Vignéras herself derives the result of \cite{Vig77} in a second paper by considering the Weil representation and mentions 
that the result is expected to hold for Hilbert and Siegel theta series as well, see \cite{Vig76}.

In the following, we prove this for the latter case by describing Siegel theta series for indefinite quadratic forms and deriving a generalization of Vignéras' result for generic genus $n$. We adopt an elementary approach similar to the one in \cite{Vig77}, which has the advantage that we explicitly construct a basis of suitable functions. This construction also embeds the known results for positive definite quadratic forms, which is for instance described by Freitag \cite{Fre83}. In this case, these "suitable functions" are harmonic polynomials and one obtains holomorphic series. However, the Siegel theta series that are constructed in the present paper are in general non-holomorphic.
In a sequel to this paper, we will investigate the special case where the quadratic form has signature $(m-1,1)$ and, by applying the result shown here, deduce the modularity of non-holomorphic Siegel theta series, which are related to holomorphic (non-modular) Siegel theta series.

We give a short overview on the main results. We use standard conventions concerning the notation, so $\e(z):=\exp(2\pi iz)$ and multiplication is hierarchically higher than division, for example $1/8\pi$ means $1/(8\pi)$.
\begin{defi}
Throughout this paper, let  $A\in \mymatrix{m}{\Z}$ denote a non-degenerate symmetric matrix of signature $(r,s)$.
\end{defi}
\begin{remark}
Note that we do not generally assume that $A$ is even. Also, in some sections we explicitly set $s=0$ and thus employ properties of the then positive definite matrix $A$.
\end{remark}
We construct modular forms on the Siegel upper half-space $$\H_n :=\lbrace Z=X+iY \mid X,Y \in \mymatrix{n}{\R}\text{ symmetric}, Y \text{ positive definite}\rbrace$$ in the form of Siegel theta series. 
We denote by $\mathcal{S}(\R^{m\times n})$ the space of Schwartz functions on $\R^{m\times n}$ and then choose $f:\R^{m\times n}\longrightarrow\R$ such that $$f(U)\ko \exp\bigl(-\pi \tr (U^\trans A U)\bigr)\in \mathcal{S}(\R^{m\times n}).$$ This ensures the absolute convergence of the theta series that we define in the following.
\begin{defi}\label{align_definition_theta}
Let $H,K\in \R^{m\times n}$ and let $\lambda\in \Z$. The theta series with characteristics $H$ and $K$ associated with $f$ and $A$ is
\begin{align*}
\mythetac{H}{K}(Z)=\theta_{H,K,f,A}(Z)
:=\det Y^{-\lambda/2}\sum\limits_{U\in H+\Z^{m\times n}}f(UY^{1/2})\ko \e\bigl( \tr(U^\trans AUZ)/2+\tr(K^\trans AU)\bigr). 
\end{align*}
\end{defi}
\begin{remark}
We drop the parameters $f$ and $A$ in the index, when the transformation of $\mythetac{H}{K}$ leaves them invariant.
In the following, it becomes clear that the choice of $\lambda$ depends on $f$, so we do not include it as additional parameter in the definition.
\end{remark}
For a positive definite matrix $A$, we consider polynomials $P:\C^{m\times n}\longrightarrow \C$ that satisfy $P(UN)=\det N^\alpha\ko P(U)$ for all $N\in \C^{n\times n}$ and a fixed $\alpha\in \N_0$. These polynomials form a complex vector space, which we denote by $\myhomopol{\alpha}$. For a modified polynomial $$p(U)=\exp\bigl(-\laplacian_A/8\pi\bigr)\bigl(P(U)\bigr)\quad\text{where}\quad \laplace_A:= \Bigl(\partialmatrix{U}\Bigr)^\trans\ko A^{-1}\partialmatrix{U},$$
and when we take $A$ to be even and set $\lambda=\alpha$, the theta series $\theta_{\mathrm{O},\mathrm{O},p,A}$ transforms like a Siegel modular form of weight $m/2+\alpha$ on a congruence subgroup of $\Gamma_n$ and with some character, where both depend on the level of $A$. If $P\in \myhomopol{\alpha}$ is annihilated by the Laplacian $\laplacian_A$, we obtain the holomorphic theta series considered by Freitag \cite{Fre83}.

When $A$ denotes an indefinite quadratic form of signature $(r,s)$, we write $A=\Apos +\Aneg$ with a positive semi-definite matrix $\Apos$ and a negative semi-definite matrix $\Aneg$ and denote by $\Amaj=\Apos -\Aneg$ the positive definite majorant matrix of $A$ (see Remark \ref{remark_decomp}). We consider the function
\begin{align*}
g(U)=\exp\bigl( -\laplacian_{\Amaj}/8\pi\bigr)\bigl(P(U)\bigr)\ko \exp\bigl(2\pi \tr(U^\trans \Aneg U )\bigr),
\end{align*}
assuming that $P\in \myhomopol{\alpha +\beta}$ factorizes as $P(U)=P_\alpha(\Upos)\cdot P_\beta(\Uneg)$ with $P_\alpha\in \myhomopol{\alpha},\ko P_\beta\in \myhomopol{\beta}$ and $U=\Upos+\Uneg$, where $\Upos$ denotes the part of $U$ that belongs to the subspace on which $A$ is positive semi-definite, i.\, e. $\tr \bigl((\Upos)^\trans A \Upos \bigr)=\tr \bigl(U^\trans \Apos U \bigr)$ and similarly $\tr \bigl((\Uneg)^\trans A \Uneg \bigr)=\tr \bigl(U^\trans \Aneg U \bigr)$.
For this choice of $g$ and considering an even matrix $A$ and setting $\lambda=\alpha-\beta-s$, the theta series $\theta_{\mathrm{O},\mathrm{O},g,A}$ transforms like a non-holomorphic Siegel modular form of weight $m/2+\lambda$ on a congruence subgroup of $\Gamma_n$ and with some character, where both depend on the level of $A$.

These explicit constructions do not only give examples of Siegel modular forms, but by applying Vignéras' result for genus $n=1$, we show that we obtain a similar criterion as in \cite{Vig77} to determine whether a Siegel theta series transforms like a modular form:
\begin{theo}\label{maintheorem}
Let $\lambda \in \Z$ and let $f:\R^{m\times n}\longrightarrow\R$ such that $$f(U)\ko \exp\bigl(-\pi \tr (U^\trans A U)\bigr)\in \mathcal{S}(\R^{m\times n})$$
and $f$ is a solution of the $n\times n$ system of partial differential equations
\begin{align*}
\Bigl(\eulerop-\frac{\laplace_A}{4\pi}\Bigr)f= \lambda\cdot I\cdot f\qquad\text{with}\quad\eulerop:= U^\trans  \frac{\partial}{\partial{U}}\quad\text{and}\quad\laplace_A\text{ as defined above}.
\end{align*}
For $H=K=\mathrm{O}$ and $A$ even, the theta series $\theta_{H,K,f,A}$ in Definition \ref{align_definition_theta} transforms like a Siegel modular form of genus $n$ and weight $m/2+\lambda$, where the level and character depend on $A$.
\end{theo}
\begin{remark}
In this paper, we determine the transformation behavior of $\theta_{H,K,f,A}$ with respect to the transformations $Z\mapsto Z+S$ for a symmetric matrix $S\in \Z^{n\times n}$ (see Lemma \ref{lemma_translation}) and $Z\mapsto -Z^{-1}$ (see Proposition \ref{proposition_transformation_indef}).
The results hold for any $H,K\in \R^{m\times n}$ and we do not generally assume that $A$ is even. By setting further preconditions for $H,K$ and $A$, one can then construct vector-valued Siegel modular forms of genus $n$ and weight $m/2+\lambda$ on the full Siegel modular group or scalar-valued modular forms on congruence subgroups, see also Remark \ref{remark_modular}. However, we will not explicitly elaborate on that here.
\end{remark}
The outline of the paper is as follows: In Section \ref{section_notation}, we briefly summarize the most important notions about Siegel modular forms that are relevant for this paper. In the next section, we examine the complex vector space formed by the solutions of the $n\times n$ system of partial differential equations from Theorem \ref{maintheorem}. Under the additional assumption that a solution $f$ must satisfy the growth condition $f(U)\ko \exp(-\pi \tr(U^\trans AU))\in \mathcal{S}(\R^{m\times n})$, we explicitly determine a basis (which is finite if $A$ is positive or negative definite and infinite otherwise) of this vector space.
In Section \ref{sectiontransformation}, we show that these basis elements can be used to construct theta series of genus $n$ that transform like Siegel modular forms of weight $m/2+\lambda$. In order to do so, we first construct non-holomorphic theta series for positive definite quadratic forms. With some modifications, this can be generalized to theta series associated with indefinite quadratic forms.
\section{Notation and preliminaries}\label{section_notation}
We fix notation and summarize standard results about Siegel modular forms and refer to Andrianov \cite[p.\,1-25]{And09} and Freitag \cite{Fre83} for further details. For convenience, we replicate some definitions of the last section. We also comment on results by Borcherds \cite{Bor98} and Vignéras \cite{Vig76} and point out differences to our set-up.

Denote the Siegel upper half-space by $$\H_n :=\lbrace Z=X+iY \mid X,Y \in \mymatrix{n}{\R}\text{ symmetric}, Y \text{ positive definite}\rbrace.$$
We let $Y^{1/2}$ denote the uniquely determined symmetric positive definite matrix that satisfies $Y^{1/2}\cdot Y^{1/2}=Y$. 
The same holds for the square root of $A$, when $A$ is a positive definite matrix.

We define modular forms on $\H_n$ for the full Siegel modular group $$\Gamma_n:=\big\lbrace M=\left(\begin{smallmatrix}
A&B\\
C&D
\end{smallmatrix}\right)\in \mymatrix{2n}{\Z} \mid M^\trans JM=J \big\rbrace,\quad\text{where}\quad J=\left(\begin{smallmatrix}
\mathrm{O}&I_n\\
-I_n&\mathrm{O}
\end{smallmatrix}\right),$$ 
which operates on $\H_n$ by
$$Z\mapsto M\langle Z\rangle=(AZ+B)(CZ+D)^{-1}.$$
The imaginary part $Y$ of $Z$ and the imaginary part $\widetilde{Y}$ of $M\langle Z\rangle$ satisfy the relation
\begin{align}\label{imaginary_part}
(C\overline{Z}+D)^\trans\ko  \widetilde{Y} (CZ+D)=Y.
\end{align}
In particular, $\widetilde{Y}$ is positive definite and symmetric.
\begin{defi}
We call $F:\H_n\longrightarrow \C$ a (classical) Siegel modular form of genus $n$ and weight $k$ if the following conditions hold:
\begin{itemize}
\item[(a)] The function $F$ is holomorphic on $\H_n$,
\item[(b)] For every $M\in \Gamma_n$ we have $F(M\langle Z\rangle)=\det (CZ+D)^{k}\ko F(Z)$,
\item[(c)] $|F(Z)|$ is bounded on domains in $\H_n$ of the form $\H^{\epsilon}_n:=\{X+iY\in \H_n\mid Y\geq \epsilon \cdot I\}$ with $\epsilon >0$.
\end{itemize}
\end{defi}
Note that the weight is not necessarily an integral number. In this context, we define -- as usual -- for $z\in \C$ and any non-integer exponent $r$ that $z^r:=\exp (r\log z)$, where $\log z=\log|z|+i\arg (z),\ko -\pi <\arg(z)\leq \pi$.
 
Due to the Koecher principle (cf. \cite[p.\,44f.]{Fre83}), which holds for $n>1$, all functions satisfying (a) and (b) admit a Fourier expansion over positive semi-definite even symmetric matrices and are in particular bounded on $\H^{\epsilon}_n$ for any $\epsilon >0$. So we do not need to impose an analogue of (c) as condition.
If we consider non-holomorphic modular forms, the Koecher principle does not necessarily hold anymore. In our case, we build Siegel theta series by using Schwartz functions and obtain absolutely convergent series, so these functions also satisfy condition (c).
\begin{remark}\label{remark_modular}
The full Siegel modular group $\Gamma_n$ is generated by the matrices
$
\bigl(\begin{smallmatrix}
I_n&S\\
\mathrm{O}&I_n
\end{smallmatrix}\bigr)$ with $S=S^\trans $ and 
$\bigl(\begin{smallmatrix}
\mathrm{O}&-I_n\\
I_n&\mathrm{O}
\end{smallmatrix}\bigr)$ (cf. \cite[p.\,322-328]{Fre83}),
so any function $F$ with $F(Z+S)=F(Z)$ for symmetric matrices $S \in \mymatrix{n}{\Z}$ and $F(-Z^{-1})=\det Z^k F(Z)$ satisfies condition (b).
For the theta series with characteristics $H,K$ that we construct here, we observe the following:
Up to a factor depending on $H$, $A$ and $S$, we can write $\mythetac{H}{K}(Z+S)$ as a theta series of the same form but with a slightly changed characteristic $H,\widetilde{K}$, see Lemma \ref{lemma_translation}. We can express $\mythetac{H}{K}(-Z^{-1})$ as a linear combination of theta series $\mythetac{J+K}{-H}(Z)$, where $J\in A^{-1}\Z^{m\times n}\operatorname{mod} \Z^{m\times n}$, see Proposition \ref{proposition_transformation_indef}.

If $A$ is an even unimodular matrix and $H=K=\mathrm{O}$, the theta series transforms like a modular form on the full group $\Gamma_n$, see Example \ref{example_positive} when $A$ is positive definite and Example \ref{example_negative} when $A$ is indefinite (in the last case we might obtain a character of $\Gamma_n$ as an additional automorphic factor).

If $H$ and $K$ are rational matrices, we can take the series  $\mythetac{H}{K}$ as entries of vector-valued functions, which then define modular forms on the full Siegel modular group. In another approach (see for example Andrianov and Maloletkin \cite{AM75}), one could consider suitable congruence subgroups of finite index in $\Gamma_n$.
\end{remark}
In Section \ref{section_vigneras} as well as Section \ref{sectiontransformation}, we will consider a fixed decomposition of the non-degenerate matrix $A$ of signature $(r,s)$, so we give a precise description here.
\begin{remark}\label{remark_decomp}
Let $\vcolumnvector{1},\ldots,\vcolumnvector{r}$ denote the eigenvectors that correspond to the positive eigenvalues of $A$ and $\vcolumnvector{r+1},\ldots,\vcolumnvector{m}$ the ones that correspond to the negative eigenvalues. We normalize these eigenvectors in a suitable way so that for
$S=(\vcolumnvector{1},\ldots,\vcolumnvector{m})\in \mymatrix{m}{\R}$
$$S^\trans AS=\I \quad\text{with}\quad\I:= \left(\begin{matrix}
I_r&\mathrm{O}\\
\mathrm{O}&-I_s
\end{matrix}\right).$$
As $\lbrace \vcolumnvector{1},\ldots, \vcolumnvector{m}\rbrace$ forms a basis of $\R^m$, we write any vector $\uvector\in\R^m$ as $\uvector=\sum_{i=1}^r \lambda_i \vcolumnvector{i}+\sum_{i=r+1}^m \lambda_i \vcolumnvector{i}$ and define $\upos:=\sum_{i=1}^r \lambda_i \vcolumnvector{i}$ and $\uneg:=\sum_{i=r+1}^m \lambda_i \vcolumnvector{i}$.

So for the inverse of $S$, we have $A=(S^{-1})^\trans\ko \I S^{-1}$.
This enables us to write $A$ as the sum of the positive semi-definite respectively negative semi-definite matrices 
\begin{align*}
\Apos:=(S^{-1})^\trans \left(\begin{matrix}
I_r&\mathrm{O}\\
\mathrm{O}&\mathrm{O}
\end{matrix}\right) S^{-1}\quad \text{and}\quad \Aneg:=(S^{-1})^\trans \left(\begin{matrix}
\mathrm{O}&\mathrm{O}\\
\mathrm{O}&-I_s
\end{matrix}\right) S^{-1}.
\end{align*}
We also associate the positive definite matrix
$\Amaj:=(S^{-1})^\trans\ko S^{-1}=\Apos -\Aneg$.
If we write $U\in \R^{m\times n}$ as $U=\Upos +\Uneg$, where $\Upos:=(\ucolumnvectorpos{1},\ldots ,\ucolumnvectorpos{n})$ and $\Uneg:=(\ucolumnvectorneg{1},\ldots ,\ucolumnvectorneg{n})$, it is straightforward to check that $$\tr \bigl((\Upos)^\trans A \Upos \bigr)=\tr \bigl(U^\trans \Apos U \bigr)\quad\text{and}\quad\tr \bigl((\Uneg)^\trans A \Uneg \bigr)=\tr \bigl(U^\trans \Aneg U \bigr).$$
\end{remark}
As our construction of Siegel theta series in Section \ref{sectiontransformation} is very similar to Borcherds' set-up \cite{Bor98} for $n=1$, we briefly recall his result and point out the main differences.
\begin{remark}\label{remark_borcherds}
Borcherds considers a non-degenerate quadratic form $Q$ with signature $(r,s)$, an even lattice $L\subset\R^m$ with the associated dual lattice $L'$ and an isometry $v$ mapping $L\otimes \R$ to $\R^{r,s}$. Considering the inverse images $v^+$ and $v^-$ of $\R^{r,0}$ and $\R^{0,s}$ under $v$, one decomposes $L\otimes \R$ in the orthogonal direct sum of a positive definite subspace $v^+$ and a negative definite subspace $v^-$. For the projection of $\boldsymbol{\lambda} \in L\otimes \R$ into $v^{\pm}$ one writes $\boldsymbol{\lambda}_{v^{\pm}}$ and obtains the positive definite quadratic form $Q_v(\boldsymbol{\lambda})=Q(\boldsymbol{\lambda}_{v^+})-Q(\boldsymbol{\lambda}_{v^-})$. As the decomposition into the subspaces $v^+$ and $v^-$ is not unique, Borcherds' theta series include an additional parameter to indicate the choice of $v^+\in G(M)$, where the Grassmannian $G(M)$ denotes the set of positive definite $r$-dimensional subspaces of  $L\otimes \R$. 
For $z \in \H_1,\hvector,\kvector\in L\otimes \R,\boldsymbol{\gamma}\in L'/L$, $\Delta$ the Laplacian on $\R^m$, and $p:\R^m\longrightarrow\R$ a polynomial that is homogeneous of degree $\alpha$ in the first $r$ variables and homogeneous of degree $\beta$ in the last $s$ variables, he defines
\begin{multline*}\theta_{L+\boldsymbol{\gamma}}(z,\hvector,\kvector;v,p):=\sum_{\boldsymbol{\lambda}\in L+\boldsymbol{\gamma}} \exp(-\Delta/8\pi y)(p)\bigl(v(\boldsymbol{\lambda}+\hvector)\bigr)\\ \cdot e\bigl(z (\boldsymbol{\lambda}+\hvector)^2_{v^+}/2+\overline{z}(\boldsymbol{\lambda}+\hvector)^2_{v^-}/2-(\boldsymbol{\lambda}+\hvector/2,\kvector)\bigr)\end{multline*}
and shows that this is a non-holomorphic modular form of weight $(r/2+\alpha,s/2+\beta)$.

In the present paper, we fix the decomposition $A=\Apos +\Aneg$ and the majorant matrix $M=\Apos - \Aneg$ by taking the eigenvectors of $A$ as a basis in $\R^m$. Then $U\in \R^{m\times n}$ is projected onto $\Upos$ in the positive definite subspace and $\Uneg$ in the negative definite subspace.
However, choosing any other decomposition of $A$ into a negative and a positive definite part leads to an analogous construction.

In Definition \ref{align_definition_theta}, we represented $\mythetac{H}{K}$ such that the analogy with Vignéras' construction (see Remark \ref{remark_vigneras}) is visible. We can also write these theta series as
\begin{multline*}
\mythetac{H}{K}(Z)=\sum_{U\in H+ \Z^{m\times n}}\exp(-\tr(\laplace_MY^{-1})/8\pi)\bigl(P(U)\bigr)\\
\cdot e\bigl(\tr(U^\trans \Apos UZ)/2+\tr(U^\trans \Aneg U\overline{Z})/2+\tr (K^\trans AU)\bigr),
\end{multline*}
which resembles Borcherds' construction. Note that we can multiply the series by $\det Y^{s/2+\beta}$ to obtain the weight $m/2+\lambda$ (where $\lambda=\alpha-\beta -s$) instead of $(r/2+\alpha,s/2+\beta)$.
\end{remark}
We conclude this section by reviewing Vignéras' construction \cite{Vig77} and addressing essential differences.
\begin{remark}\label{remark_vigneras}
Vignéras considers theta series of genus 1
\begin{align*}
\mythetac{\boldsymbol{0}}{\boldsymbol{0}}(z)=y^{-\lambda/2}\sum\limits_{\uvector\in L}f(\uvector\sqrt{y})\ko \e\bigl( Q(\uvector)z\bigr),
\end{align*}
where $L\subset \R^m$ denotes a lattice, $Q(\uvector)=\frac12\ko \uvector^\trans A \uvector$ a quadratic form of signature $(r,s)$ and $z=x+iy$ an element of the upper half-plane $\H_1$.
The following two requirements are imposed on the function $f$: Set
$\widetilde{f}(\uvector)=f(\uvector)\ko  \exp\bigl(-2\pi Q(\uvector)\bigr)$. Then
for any polynomial $p:\R^{m}\longrightarrow\R$ with $\deg(p)\leq 2$ and any partial derivative $\partial^{\alpha}$ with $|\alpha|\leq 2$,
\[p\cdot\widetilde{f} \in \mathcal{L}^1(\R^{m})\cap\mathcal{L}^2(\R^{m})\quad\text{and}\quad \partial^{\alpha}\widetilde{f} \in \mathcal{L}^1(\R^{m})\cap\mathcal{L}^2(\R^{m}).\]
Furthermore, $f$ satisfies the differential equation of second order
$$\Bigl(E-\frac{\Delta_A}{4\pi}\Bigr)f=\lambda\cdot  f\qquad\text{with}\quad E:=\sum_{d=1}^m u_d\ko  \frac{\partial}{\partial {u_d}}\quad\text{and}\quad \Delta_A:=\sum\limits_{a,b=1}^m\partialmatrix{u_{a}} (A^{-1})_{ab} \partialmatrix{u_{b}}.$$
Then $\mythetac{\boldsymbol{0}}{\boldsymbol{0}}$ transforms like a modular form of weight $m/2+\lambda$.
\end{remark}
For higher genus $n\in \N$, we introduce some notation to formulate an analogous growth condition.
For $p\in [1,\infty)$ let $\mathcal{L}^p(\R^{m\times n})$ denote the Lebesgue space of functions $f:\R^{m\times n}\longrightarrow \C$ for which 
\begin{align*}
\|f\|_p:=\Biggl(\ko \myintegral |f(U)|^p\ko dU\Biggr)^{1/p}
\end{align*}
is finite.
We use the usual multi-index notation on $\R^{m\times n}$, where
$\alpha\in \N_0^{m\times n}$ with $|\alpha|=\sum_{i=1}^{m}\sum_{j=1}^{n} \alpha_{ij}$, so
$$U^{\alpha}=\prod\limits_{\substack{1\leq i\leq m\\ 1\leq j\leq n}} U_{ij}^{\alpha_{ij}}\quad\text{and}\quad
\partial^{\alpha}=\prod\limits_{\substack{1\leq i\leq m\\ 1\leq j\leq n}}\Bigl(\frac{\partial}{\partial U_{ij}}\Bigr)^{\alpha_{ij}}.$$
For $f:\R^{m\times n}\longrightarrow \R$, one sets $\widetilde{f} (U):=f(U)\ko \exp\bigl(-\pi \tr (U^\trans A U)\bigr)$ and -- analogously to Vignéras -- assumes that for any polynomial $p:\R^{m\times n}\longrightarrow\R$ with $\deg(p)\leq 2$ and any partial derivative $\partial^{\alpha}$ with $|\alpha|\leq 2$,
\begin{align}\label{align_convergence_assumption}
p\cdot\widetilde{f} \in \mathcal{L}^1(\R^{m\times n})\cap\mathcal{L}^2(\R^{m\times n})\quad\text{and}\quad \partial^{\alpha}\widetilde{f} \in \mathcal{L}^1(\R^{m\times n})\cap\mathcal{L}^2(\R^{m\times n}).
\end{align}
This allows us to apply Vignéras' result for theta series of genus 1 (as we make use of the fact that Hermite functions build an orthogonal basis of $\mathcal{L}^2$-functions) and the Poisson summation formula.

However, for simplification, we replace assumption \eqref{align_convergence_assumption} by the more restrictive assumption that $\widetilde{f}$ is a Schwartz function.
\section{A generalization of Vignéras' differential equation}\label{section_vigneras}
To derive an analogue of Vignéras' result for Siegel modular forms of higher genus $n\in \N$, we introduce matrix-valued operators generalizing $E$ and $\Delta_A$. 
\begin{defi}
For $U\in \R^{m\times n}$ let $\partial/\partial U=\bigl(\partial/\partial U_{\mu \nu}\bigr)_{1\leq \mu\leq m,\ko1\leq \nu\leq n}$.
We define the generalized Euler operator
\[\eulerop:= U^\trans  \frac{\partial}{\partial{U}}\quad\text{with}\quad\eulerentry=\sum\limits_{d=1}^m U_{di}\ko  \frac{\partial}{\partial U_{dj}}\quad (1\leq i\leq n,\,1\leq j\leq n)\]
and the generalized Laplace operator associated with $A$
\begin{align*}
\laplace_A := \Bigl(\partialmatrix{U}\Bigr)^\trans\ko A^{-1}\ko\partialmatrix{U}\quad\text{with}\quad (\laplace_A)_{ij}=\sum\limits_{a,b=1}^m \frac{\partial}{\partial U_{ai}}(A^{-1})_{ab}\frac{\partial}{\partial U_{bj}}.
\end{align*}
For the normalized Laplacian $\laplace_I$ we simply write $\laplace$.
Further, we set 
\begin{align*}
    \mathcal{D}_A:=\eulerop-\frac{\laplace_A}{4\pi}.
\end{align*}
\end{defi}
The $n\times n$ system of partial differential equations
\begin{align}\label{align_dgl_A}
\mathcal{D}_A f= \lambda\cdot I\cdot f\quad\text{for }\lambda\in\Z\text{ and }A \text{ indefinite of signature }(r,s)
\end{align}
is a direct generalization of the set-up in \cite{Vig77}. In this section, we examine the complex vector space formed by the solutions $f$ of \eqref{align_dgl_A} that additionally satisfy the growth condition $$f(U)\ko \exp(-\pi \tr(U^\trans AU))\in \mathcal{S}(\R^{m\times n}).$$ We explicitly determine a basis (which is finite if $A$ is positive or negative definite and infinite otherwise) of this vector space.
\subsection{Functions with a homogeneity property}\label{section_homogeneity}
As mentioned in the introduction, we employ polynomials with a certain homogeneity property to construct Siegel theta series. In the following, we introduce the complex vector space of all functions with this homogeneity property. For a \emph{differentiable} function $f$, we show in Proposition \ref{proposition_homogen_dgl} that $f$ is homogeneous of degree $\alpha$ if and only if $f$ solves the system of partial differential equations $\eulerop f=\alpha\cdot I\cdot f$. Further, we show in Lemma \ref{lemma_konstante} that for a \emph{polynomial} function $p$ it is already sufficient that $\eulerop p=C\cdot p$ holds for some $C\in \C^{n\times n}$ to deduce that $p$ is a homogeneous function.
\begin{defi}
For $\alpha\in \N_0$, $m,n\in \N$, we define the complex vector space
\begin{align*}
\myhomofct{\alpha}:=\lbrace f: \C^{m\times n}\longrightarrow \C \mid f \text{ continuous},\, f(UN)=\det N^\alpha f(U)\text{ for all }N \in \mymatrix{n}{\C}\rbrace.
\end{align*}
For $n=1$, this is the usual definition of a homogeneous function of non-negative degree.  
As a subspace, we consider all polynomials of this class, which is the space $\myhomopol{\alpha}$ from the introduction.
\end{defi}
\begin{remark}
The vector space $\myhomopol{\alpha}$ is described by Maass \cite{Maa59}. He determines the structure of $\myhomopol{\alpha}$, shows that it has finite dimension and even gives an explicit formula for the dimension. In the following, $\basis{\alpha}{n}$ denotes a finite basis of $\myhomopol{\alpha}$. We state some observations to show that we obtain non-trivial examples.
\begin{itemize}
\item For $m<n$, we have $\myhomofct{\alpha}= \C$: We take $U\in \R^{m\times n}$ such that $f(U)\neq 0$. One can multiply elementary matrices from the right such that $U$ is in reduced column echelon form. If $U$ has less rows than columns, at least the last column is a zero column. Setting $N=\operatorname{diag}(1,\ldots,1,\lambda)$ with $\lambda\notin\lbrace 0,1\rbrace$, leads to the identity $f(U)=\lambda^\alpha f(U)$, 
which is only satisfied for $\alpha=0$. The orbit of the right action of invertible matrices on $U\in\C^{m\times n}$ is dense and $f$ is continuous, so $f$ is a constant function.
\item Note that $f\cdot g\in \myhomofct{\alpha+\beta}$ for $f\in \myhomofct{\alpha},\ko  g\in \myhomofct{\beta}$ and  $f+ g\in \myhomofct{\alpha}$ for $f,g\in \myhomofct{\alpha}$.
\item For $m\geq n$ let $\widetilde{U}\in \mymatrix{n}{\C}$ be a square submatrix of maximal size of $U\in \C^{m\times n}$. Clearly, we have $\det \widetilde{U}^{\alpha}\in\myhomopol{\alpha}$. Due to this and by picking up on the previous point, we obtain all functions in $\myhomopol{\alpha}$ by taking the product of $\alpha$ (possibly different) $n\times n$--minors $\det \widetilde{U}$ and linear combinations thereof.
\end{itemize}
\end{remark}
Homogeneous functions that are also differentiable are characterized by the identity $Ef=\alpha\cdot f$. We observe that this statement can be generalized.
\begin{satz}\label{proposition_homogen_dgl}
Let $f: \C^{m\times n}\longrightarrow \C$ be a differentiable function. We have $f\in \myhomofct{\alpha}$ if and only if $\eulerop f=\alpha\cdot I\cdot f$.
\end{satz}
\begin{proof}
For $U\in \C^{m\times n}$ and $N\in \C^{n\times n}$, the derivative of the entry $(UN)_{k\ell}=\sum\limits_{\nu =1}^n U_{k\nu}N_{\nu \ell}$ with $1\leq k\leq m,\,1\leq \ell \leq n$ is $$\frac{\partial (UN)_{k\ell}}{\partial N_{ij}}=\begin{cases}
0&\text{ if }j \neq \ell,\\
U_{ki}&\text{ if }j=\ell.
\end{cases}$$
Therefore,
\begin{align*}
\partialmatrix{N_{ij}} \bigl(f(UN)\bigr)=\sum\limits_{k=1}^m \sum\limits_{\ell=1}^n \frac{\partial f}{\partial U_{k\ell}}(UN)\ko  \frac{\partial (UN)_{k\ell}}{\partial N_{ij}}=\sum\limits_{k=1}^m U_{ki}\ko  \frac{\partial f}{\partial U_{kj}} (UN).
\end{align*}
Hence, we obtain for the derivative of $f(UN)$ with respect to $N$ that $$\partialmatrix{N} \bigl(f(UN)\bigr)=U^\trans \ko  \frac{\partial f}{\partial U} (UN)$$ and by the definition of the generalized Euler operator $\eulerop$ with respect to $U$ it is
\begin{align}\label{aligneuler}
(\eulerop f)(UN)=(UN)^\trans \ko  \frac{\partial f}{\partial U}(UN)=N^\trans \ko  \partialmatrix{N} \bigl(f(UN)\bigr).
\end{align}
The adjugate matrix $\operatorname{adj}(N)\in \mymatrix{n}{\C}$ is defined as $\bigl(\operatorname{adj}(N)\bigr)_{ij}:=(-1)^{i+j} \det {\widetilde{N}}_{ji},$ where ${\widetilde{N}}_{ji}$ denotes the $(n-1)\times (n-1)$-matrix obtained by deleting the $j$-th row and $i$-th column. Laplace expansion of the determinant gives
\begin{align*}
\partialmatrix{N_{ij}} (\det N)=\partialmatrix{N_{ij}} \biggl( \sum\limits_{k=1}^n(-1)^{i+k}N_{ik}\det {\widetilde{N}}_{ik} \biggr)=(-1)^{i+j}\det {\widetilde{N}}_{ij}.
\end{align*}
Hence, the derivate of the determinant of $N$ is the transpose of the adjugate matrix:
\begin{align}\label{aligndeterminant}
\partialmatrix{N} (\det N) = \operatorname{adj}(N)^\trans
\end{align}
For $f\in \myhomofct{\alpha}$ the identity $f(UN)=\det N^\alpha f(U)$ holds for all $N\in \mymatrix{n}{\C}$. From equations (\ref{aligneuler}) and (\ref{aligndeterminant}) it follows
\begin{multline*}
\bigl(\eulerop f\bigr)(UN)= N^\trans  \partialmatrix{N} \bigl(f(UN)\bigr)=N^\trans  \partialmatrix{N}\bigl(\det N^\alpha f(U)\bigr)\\
=\alpha \det N^{\alpha-1} \cdot N^\trans    \operatorname{adj}(N)^\trans \cdot f(U)=\alpha \det N^\alpha \cdot I\cdot f(U),
\end{multline*}
since $\operatorname{adj}(N) \ko  N=\det N\cdot I$. We set $N=I$ and obtain the identity $\eulerop f=\alpha\cdot I\cdot f$.\\
To show the other implication, notice that $f(UN)(\det N)^{-\alpha}$ is constant with respect to $N$ if $f$ satisfies $\eulerop f=\alpha\cdot I\cdot f$: Using the identities (\ref{aligneuler}) and (\ref{aligndeterminant}), we obtain
\begin{align*}
N^\trans  \partialmatrix{N} \bigl( f(UN) \det N^{-\alpha}\bigr)\vspace{2pt}
=\det N^{-\alpha} \cdot\bigl(\eulerop f\bigr)(UN)-\alpha \det N^{-\alpha}\cdot I\cdot f(UN)=\mathrm{O}.
\end{align*}
Thus, we have $f(UN) \det N^{-\alpha}=C(U),$ where $C$ is independent of $N$. For $N=I$, this is $f(U)$, and hence we conclude $f(UN)=\det N^\alpha f(U)$.\qedhere
\end{proof}
Later we will only consider polynomial solutions. In this case, we can state the following lemma, which can be left aside for now, but will be used in the proof of Proposition \ref{corollary_indef}.
\begin{lemma}\label{lemma_konstante}
Let $p:\C^{m\times n}\longrightarrow\C$ be a polynomial that solves the system of partial differential equations
\begin{align}\label{align_dgl_konstante}
\eulerop p=C\cdot p\quad(C\in \C^{n\times n}).
\end{align}
If $p$ is not the zero function, the matrix $C$ has the form $C=\alpha\cdot I$ for some $\alpha\in\N_0$.
\end{lemma}
\begin{proof}
First, we examine the case $m=n=2$ and write $U=\left(\begin{smallmatrix} a&b\\ c&d
\end{smallmatrix}\right)$ and
\begin{align*}
p(a,b,c,d)=\sum\limits_{\alpha,\beta,\gamma,\delta \in \N_0} c_{\alpha,\beta,\gamma,\delta}\ko a^{\alpha} b^\beta c^\gamma d^\delta\quad\text{with}\quad c_{\alpha,\beta,\gamma,\delta}\in \C.
\end{align*}
By assumption, $p$ satisfies the $2\times 2$-system of partial differential equations
\begin{align}\label{align_22system}
\begin{pmatrix}
a&c\\
b&d
\end{pmatrix}
\begin{pmatrix}
\partialmatrix{a}&\partialmatrix{b}\\[2pt]
\partialmatrix{c}&\partialmatrix{d}
\end{pmatrix}p=\begin{pmatrix}
C_{11}&C_{12}\\
C_{21}&C_{22}
\end{pmatrix}\cdot p.
\end{align}
Considering the upper left equation, we have
\begin{align*}
\Bigl(a\partialmatrix{a}+c\partialmatrix{c}\Bigr)p= \sum\limits_{\alpha,\beta,\gamma,\delta \in \N_0} (\alpha+\gamma)\ko c_{\alpha,\beta,\gamma,\delta}\ko a^{\alpha} b^\beta c^\gamma d^\delta =C_{11}\cdot p,
\end{align*}
thus $\alpha+\gamma=C_{11}$. Analogously, we deduce by the bottom right equation that $\beta+\delta=C_{22}$ holds. As $p$ is a polynomial, $C_{11}$ and $C_{22}$ denote non-negative integers. We write $C_{11}=k$ and $C_{22}=\ell$ from now on. We have shown that $p$ is homogeneous (in the original sense) of degree $k$ in the variables of the first column $a,c$ and homogeneous of degree $\ell$ in the variables of the last column $b,d$. It is easy to see that $C_{12}=C_{21}=0$ holds:
By assumption, the polynomial $p$ satisfies the upper right equation
\begin{align}\label{align_upper_right}
\Bigl(a\partialmatrix{b}+c\partialmatrix{d}\Bigr)p=C_{12}\cdot p.
\end{align}
As the left-hand side is a polynomial, homogeneous of degree $k+1$ in $a,c$ and of degree $\ell-1$ in $b,d$, and the right-hand side is a multiple of $p$, i.\,e. homogeneous of degree $k$ and $\ell$, we deduce that $C_{12}$ must equal zero. Analogously, we conclude by the bottom left equation of \eqref{align_22system} that $C_{21}=0$.

It remains to be shown that $k=\ell$ holds. We write $$p(a,b,c,d)=\sum_{\alpha+\gamma=k} a^\alpha c^\gamma\ko p_{\alpha,\gamma}(b,d),$$ where $p_{\alpha,\gamma}$ denote homogeneous polynomials in $b,d$ of degree $\ell$. Then equation (\ref{align_upper_right}) with $C_{12}=0$ has the form
\begin{align*}
\sum\limits_{\alpha +\gamma =k} a^{\alpha +1}c^\gamma\ko  \partialmatrix{b}p_{\alpha,\gamma}(b,d)+\sum\limits_{\alpha +\gamma =k} a^{\alpha }c^{\gamma+1}\ko  \partialmatrix{d}p_{\alpha,\gamma}(b,d)\equiv 0.
\end{align*}
We obtain by comparison of the coefficients of $a^\nu c^\mu,\, 0\leq \nu\leq k+1,\,\mu =k+1-\nu$:
\begin{align*}
\partialmatrix{d}p_{0,k}(b,d)\equiv 0\qquad\text{and}\qquad
\partialmatrix{d}p_{\alpha,\gamma}(b,d)=-\partialmatrix{b}p_{\alpha-1,\gamma+1}(b,d)\quad
\text{for}\quad 1\leq \alpha \leq k,\,\gamma=k-\alpha
\end{align*}
Thus, we recursively determine the structure of $p_{\alpha,\gamma}$ to be
$p_{\alpha,\gamma}(b,d)=\sum_{r=0}^\alpha e_r\ko b^{\ell-r}d^r$ with $e_r\in \R$. In particular, we see that the exponent of $d$ does not exceed $\alpha$, i.\,e. $\delta\leq \alpha$.

We make use of the symmetric structure of the polynomial $p$ and exchange $a$ and $c$ and also $b$ and $d$ in the equations above. Then we obtain $\beta\leq \gamma$.
By interchanging $a$ and $d$ along with their exponents $\alpha$ and $\delta$ as well as $b$ and $c$ along with their exponents $\beta$ and $\gamma$ and using the bottom left equation of \eqref{align_22system}, we obtain $\alpha\leq \delta$ and $\gamma\leq \beta$.
We have shown $\alpha=\delta$ and $\gamma=\beta$, and in particular, $k=\ell$ holds.

For generic $m,n\in\N$, we reduce the $n\times n$-system $\eulerop p=C\cdot p$ to the case $m=n=2$.
We write $U=(\ucolumnvector{1},\ldots,\ucolumnvector{n})$ with $\ucolumnvector{i}\in\C^m$ and choose $N\in \C^{n\times n}$ such that the $i$-th column of $U$ is substituted by $a\ucolumnvector{i}+c\ucolumnvector{j}$ and the $j$-th column by $b\ucolumnvector{i}+d\ucolumnvector{j}$, where we assume that $i<j$, i.\,e. we have
$$UN=(\ucolumnvector{1},\ldots,a\ucolumnvector{i}+c\ucolumnvector{j},\ldots,b\ucolumnvector{i}+d\ucolumnvector{j},\ldots,\ucolumnvector{n}).$$
A simple calculation yields
\begin{align*}
\begin{pmatrix}
a\partialmatrix{a}+c\partialmatrix{c}&a\partialmatrix{b}+c\partialmatrix{d}\\[2pt]
b\partialmatrix{a}+d\partialmatrix{c}&b\partialmatrix{b}+d\partialmatrix{d}
\end{pmatrix}
\bigl(p(UN)\bigr)=
\left(\begin{pmatrix}
\eulerop_{ii}&\eulerop_{ij}\\
\eulerop_{ji}&\eulerop_{jj}
\end{pmatrix} p\right)(UN).
\end{align*}
As $p$ solves \eqref{align_dgl_konstante} by assumption, we have 
\begin{align*}
\begin{pmatrix}
a&c\\
b&d
\end{pmatrix}
\begin{pmatrix}
\partialmatrix{a}&\partialmatrix{b}\\[2pt]
\partialmatrix{c}&\partialmatrix{d}
\end{pmatrix}
\bigl(p(UN)\bigr)=
\begin{pmatrix}
C_{ii}&C_{ij}\\
C_{ji}&C_{jj}
\end{pmatrix}\cdot p(UN).
\end{align*}
For $2\times 2$-systems of this form we have shown above that $C_{ii}=C_{jj}=\alpha$ for some $\alpha\in \N_0$ and  $C_{ij}=C_{ji}=0$. As we can choose any $i,j\in\lbrace 1,\ldots, n\rbrace$ with $i< j$, we deduce the claim.
\end{proof}
\subsection{Description of theta series with modular transformation behavior by partial differential equations}\label{section_coeff_diff_eq_pos}
In this section, we show the connection between the functions with the homogeneity property that was described in the last section and the functions that are employed in Section \ref{sectiontransformation} to construct modular Siegel theta series. Moreover, we apply Vignéras' result for $n=1$ to explicitly give a basis for the vector space of solutions of \eqref{align_dgl_A} under the additional growth condition.

First, we state a lemma that holds for any symmetric non-degenerate matrix $A$ of signature $(r,s)$. 
Namely, we compute the commutator of the $k$-th power of the Laplacian $(\laplacian_A)^k$ (we will drop the brackets and write $\laplacian_A^k$ for simplicity) and the Euler operator.
\begin{lemma}\label{lemma_commutator}
The commutator of $\eulerentry\,(1\leq i \leq n,\,1\leq j\leq n)$ and $\laplacian_A^k$ is 
\begin{align*}
\bigl[\eulerentry,\laplacian_A^k\bigr]:=\eulerentry \cdot \laplacian_A^k - \laplacian_A^k \cdot \eulerentry = - 2 k\cdot (\laplace_A)_{ij} \cdot\laplacian_A^{k-1}\quad(k\in\N).
\end{align*}
\end{lemma}
\begin{proof}
We show the claim by induction on $k$. For $k=1$ one calculates the commutator of $\laplacian_A$ and $\eulerentry$. 
By definition we have
\[\laplacian_A \cdot \eulerentry = \sum\limits_{c=1}^n\sum\limits_{a,b,d=1}^m \partialmatrix{U_{ac}} (A^{-1})_{ab} \partialmatrix{U_{bc}} U_{di} \partialmatrix{U_{dj}},\]
which we can write -- denoting by $\delta_{ij}$ the Kronecker delta -- as
\begin{align*}
\sum\limits_{c=1}^n &\sum\limits_{a,b,d=1}^m (A^{-1})_{ab}\ko \left( U_{di} \frac{\partial^3}{\partial U_{ac}\partial U_{bc} \partial U_{dj}}+ \delta_{ad} \delta_{ci} \frac{\partial^2}{\partial U_{bc} \partial U_{dj}}+ \delta_{bd} \delta_{ci} \frac{\partial^2}{\partial U_{ac} \partial U_{dj}} \right)\\
&\qquad=
\sum\limits_{c=1}^n\sum\limits_{a,b,d=1}^m U_{di} \partialmatrix{U_{dj}}(A^{-1})_{ab} \frac{\partial^2}{\partial U_{ac}\partial U_{bc}}
+\sum\limits_{a,b=1}^m (A^{-1})_{ab}\ko\left( \frac{\partial^2}{\partial U_{aj} \partial U_{bi}} + \frac{\partial^2}{\partial U_{bj} \partial U_{ai}}\right).
\end{align*}
Since $A^{-1}$ is symmetric, this is $\eulerentry \cdot \laplacian_A + 2 \cdot(\laplace_A)_{ij}$. 
The operators $(\laplace_A)_{ij}$ and $\laplacian_A$ commute; thus we deduce for $k\mapsto k+1$
\begin{align*}
\laplacian_A^{k+1}\cdot  \eulerentry &= \laplacian_A\cdot  \bigl(\eulerentry\cdot \laplacian_A^k+ 2 k\cdot (\laplace_A)_{ij} \cdot\laplacian_A^{k-1} \bigr)\vspace{2pt}\\
&=\bigl( \eulerentry\cdot \laplacian_A+2\cdot (\laplace_A)_{ij} \bigr)\cdot \laplacian_A^k+2k\cdot (\laplace_A)_{ij}\cdot \laplacian_A^k\vspace{2pt}\\
&=\eulerentry\cdot \laplacian_A^{k+1}+2 (k+1)\cdot(\laplace_A)_{ij}\cdot \laplacian_A^k.\qedhere
\end{align*}
\end{proof}
We can now conclude that all solutions of (\ref{align_dgl_A}) can be ascribed to functions that have the homogeneity property of degree $\lambda$ by applying the previous lemma and Proposition \ref{proposition_homogen_dgl}. 
\begin{lemma}\label{differentialequation_fr}
Let $f,g:\R^{m\times n}\longrightarrow \R$ denote functions for which $\exp (c_1\laplacian_A )f$ and $\exp ( c_2\laplacian_A)g$ are well-defined for any $c_1,c_2\in\R$ (we apply this result to polynomials $f$ and $g$ in the following, hence these conditions make sense). Moreover, we assume that $f$ and $g$ are related by $f=\exp (-\operatorname{tr}\laplace_A/8\pi)g$.
Then $f$ is a solution of \eqref{align_dgl_A} if and only if $g$ satisfies $\eulerop g=\lambda\cdot I\cdot g$, i.\,e. $g\in \myhomofct{\lambda}$.
\end{lemma}
\begin{proof}
We set $c:=-1/8\pi$ to shorten notation.
For
\begin{align*}
f(U)=\exp (c\operatorname{tr}\laplace_A)\bigl(g(U)\bigr)=\sum\limits_{k=0}^\infty \frac{c^k}{k!}\laplacian_A^k  \bigl(g(U)\bigr)\quad(U\in \R^{m\times n})
\end{align*}
we consider the entry $(i,j)$ for $1\leq i \leq n,\ko 1\leq j\leq n$ of the system of partial differential equations \eqref{align_dgl_A}:
\begin{align*}
\eulerentry f+2c\cdot(\laplace_A)_{ij} f&=\sum\limits_{k=0}^\infty \frac{c^k}{k!}\ko \bigl(\eulerentry\cdot  \laplacian_A^k \bigr) g+2\sum\limits_{k=0}^\infty \frac{c^{k+1}}{k!}\ko\bigl( (\laplace_A)_{ij} \cdot \laplacian_A^k \bigl) g \vspace{2pt}\\
&=\eulerentry g+\sum\limits_{k=1}^\infty \frac{c^k}{k!}\ko \bigl(\eulerentry\cdot  \laplacian_A^k \bigr) g+2\sum\limits_{k=1}^\infty \frac{c^k}{(k-1)!}\ko\bigl( (\laplace_A)_{ij} \cdot \laplacian_A^{k-1} \bigl) g \vspace{2pt}\\
&= \eulerentry g+\sum\limits_{k=1}^\infty \frac{c^k}{k!}\ko\bigl(\eulerentry \cdot \laplacian_A^k+2k\cdot (\laplace_A)_{ij}\cdot \laplacian_A^{k-1} \bigr)g
\end{align*}
Due to Lemma \ref{lemma_commutator}, we have
\begin{align*}
\bigl(\eulerentry \cdot \laplacian_A^k+2k\cdot (\laplace_A)_{ij}\cdot \laplacian_A^{k-1} \bigr)g=\bigl(\laplacian_A^k\cdot \eulerentry \bigr)g
\end{align*}
and therefore obtain
\begin{align*}
\eulerentry f+2c\cdot (\laplace_A)_{ij} f=\sum\limits_{k=0}^\infty \frac{c^k}{k!}\ko\laplacian_A^k\bigl( \eulerentry g\bigr)=\exp \bigl(c\laplacian_A\bigr)(\eulerentry g).
\end{align*}
If $\eulerentry g=\lambda\cdot \delta_{ij}\cdot g$ holds, the right-hand side equals $\lambda\cdot \delta_{ij}\cdot f$. As we have $g=\exp (-c\laplacian_A)f$ (see Property \ref{propertylaplace1}), we deduce that $\exp(c\laplacian_A)(\eulerentry g)=\lambda \cdot \delta_{ij}\cdot f$ implies $\eulerentry g=\lambda\cdot \delta_{ij}\cdot g$. By Proposition \ref{proposition_homogen_dgl}, this is equivalent to $g\in \myhomofct{\lambda}$.
\end{proof}
In the next proposition, we consider \eqref{align_dgl_A} for \emph{positive definite} matrices $A$, namely the system of partial differential equations
\begin{align}\label{dgl+}
\mathcal{D}_A\ko p=\alpha\cdot I\cdot p\quad\text{for $\alpha\in \N_0$ and $A$ positive definite.}
\end{align}
We determine a finite basis of all solutions of \eqref{dgl+} by additionally imposing a certain growth condition. Together with Proposition \ref{proposition_transformation_posdef}, where it is shown that theta series associated with these functions transform like modular forms, we obtain Theorem \ref{maintheorem} for positive definite matrices $A$.
In the proof, we employ Vignéras' result \cite{Vig77} and the fact that we can explicitly construct a finite basis $\basis{\alpha}{n}$ of $\myhomopol{\alpha}$ due to Maass' result \cite{Maa59}.
\begin{satz}\label{lemma_equivalence}
Let $\basis{\alpha}{n}$ denote a finite basis of $\myhomopol{\alpha}$ and let $A\in\Z^{m\times m}$ denote a positive definite symmetric matrix. 
Every solution $f:\R^{m\times n}\longrightarrow\R$ of \eqref{dgl+} that additionally satisfies the growth condition $\widetilde{f}(U):=f(U)\ko \exp \bigl( -\pi \tr (U^\trans A U) \bigr)\in \mathcal{S}(\R^{m\times n})$ is a polynomial. Moreover, a finite basis of this space of solutions is given by
\[p=\exp \bigl(-\operatorname{tr}\laplace_A/8\pi\bigr) P\quad\text{with } P\in \basis{\alpha}{n}.\]
\end{satz}
\begin{proof}
We give a short review of Vignéras' reasoning and apply it to functions with matrix variables. We identify $\R^{m\times n}$ with $\R^{m n}$ by writing $\R^{m\times n}\ni U=(\ucolumnvector{1},\ldots,\ucolumnvector{n}),$ $\ucolumnvector{i}\in \R^m,$ as column vector 
\begin{align*}
\boldsymbol{u}=\begin{pmatrix}
\ucolumnvector{1}\\
\vdots\\
\ucolumnvector{n}
\end{pmatrix}\in \R^{mn}.
\end{align*}
If $f$ satisfies the system of differential equations (\ref{dgl+}) for $c=-1/8\pi$, it follows in particular that 
\begin{align}\label{align_differentialequation_mn}
\tr (\mathcal{D}_A)f=\alpha n\cdot f
\end{align}
holds.
We have $$\tr \eulerop =\sum_{i=1}^{n}\sum_{d=1}^{m}U_{di}\ko \partialmatrix{U_{di}}\quad\text{and}\quad \tr \laplace_A =\sum\limits_{\nu=1}^n \sum\limits_{\mu,\rho =1}^m \frac{\partial}{\partial U_{\mu\nu}}\ko (A^{-1})_{\mu \rho}\ko \frac{\partial}{\partial U_{\rho\nu}},$$ which are the usual Euler operator on $\R^{m n}$ and the Laplacian associated with the positive definite $mn\times mn$- matrix that consists of blocks of $m\times m$-matrices that are zero except for $n$ copies of $A$ on the diagonal.
We write $U$ in a suitable basis such that we consider the quadratic form $(S^{-1})^\trans\ko A S^{-1}=I$ to express
 $\widetilde{f}$ in an orthogonal basis of Hermite functions $H_{\boldsymbol{k}}$ in $m n$ variables as $$\widetilde{f}=\sum_{\boldsymbol{k}\in \N_0^{mn}} c_{\boldsymbol{k}}\ko  H_{\boldsymbol{k}}\quad\text{with}\quad c_{\boldsymbol{k}}\in \R\text{ and } \boldsymbol{k}=(k_{\mu\nu})_{1\leq \mu \leq m, 1\leq \nu \leq n},$$ 
where the Hermite functions in several variables are defined in terms of Hermite functions in one dimension:
\begin{align*}
H_{\boldsymbol{k}}(U)=\prod_{\mu=1}^{m}\prod_{\nu=1}^{n}H_{k_{\mu\nu}}(U_{\mu\nu})\quad\text{with }H_{k_{\mu\nu}}(U_{\mu\nu})=\exp \bigl(\pi U_{\mu\nu}^2\bigr)\frac{d^{k_{\mu\nu}}}{dU_{\mu\nu}^{k_{\mu\nu}}} \exp\bigl(-2\pi U_{\mu\nu}^2\bigr)
\end{align*}
Since $f$ is a solution of (\ref{align_differentialequation_mn}), a basis of all functions $\widetilde{f}$ is determined by the finite set of Hermite functions $H_{\boldsymbol{k}}$ with $|\boldsymbol{k}|=\sum_{\mu=1}^{m}\sum_{\nu=1}^{n} k_{\mu\nu}=\alpha n$ (this is Vignéras' argument, see \cite{Vig77}), where we can rewrite $H_{\boldsymbol{k}}(U)=p_{\boldsymbol{k}}(U)\ko \exp\bigl(-\pi\tr(U^\trans AU)\bigr)$ with the Hermite polynomial $p_{\boldsymbol{k}}$. So $f(U)=\widetilde{f}(U)\ko \exp\bigl(\pi\tr(U^\trans A U)\bigr)$ can be expanded in terms of finitely many orthogonal Hermite polynomials $p_{\boldsymbol{k}}$ and thus is a polynomial itself.

Thus, $g:=\exp\bigl(\laplacian_A/8\pi\bigr)f$ is a polynomial that satisfies $\eulerop g=\alpha\cdot I\cdot g$ by Lemma \ref{differentialequation_fr}. We can choose any basis $\basis{\alpha}{n}$ of $\myhomopol{\alpha}$ to describe these homogeneous polynomials. Hence, we also obtain a basis of the solutions of (\ref{dgl+}). As $\myhomopol{\alpha}$ is a finite-dimensional vector space, the basis $\basis{\alpha}{n}$ is finite.
\end{proof}
Now we let $A$ denote an indefinite matrix of signature $(r,s)$ again.
When we consider the associated system of partial differential equations \eqref{align_dgl_A}, the solutions, which we describe in Proposition \ref{corollary_indef}, can be traced to functions that are defined on $U^\pm$ respectively, where $\Upos$ is the projection of $U$ onto the subspace where $A$ is positive definite and $\Uneg$ the projection into the subspace where $A$ is negative definite.
So we first consider $\I=\bigl(\begin{smallmatrix}I_r&\mathrm{O}\\
\mathrm{O}&-I_s\end{smallmatrix}\bigr)$ instead of $A$ and
the corresponding system of partial differential equations
\begin{align}\label{align_dgl}
\mathcal{D}_{\I}f= \lambda\cdot I\cdot f\quad(\lambda\in\Z),
\end{align} 
which can easily be split up into one part that depends on the first $r$ rows of $U$ and another part depending on the last $s$ rows of $U$. We write $\ur$ and $\us$ for these projections of $U$. Here, we have $M=I$ and thus $\laplace_M=\laplace$, and we show that a basis of all solutions is given by the functions
\begin{align}\label{align_solutions}
\exp\bigl( -\laplacian/8\pi\bigr)\bigl(P(U)\bigr)\ko \exp\bigl(-2\pi \tr(\us^\trans \us)\bigr),
\end{align} 
where $P$ splits as $P(U)=P_r(\ur)\cdot P_s(\us)$ with $P_r\in \basis{\alpha}{n}\subset \myhomopol{\alpha}$ and $P_s\in \basis{\beta}{n}\subset \myhomopol{\beta}$ with $(\alpha,\beta)\in \N_0^2$ such that $\alpha-\beta=\lambda+s$.
\begin{lemma}\label{lemma_produktregel}
If one applies the Laplacian $\laplace_\I$ and the Euler operator on a product of functions $g,h:\R^{m\times n}\longrightarrow \R$, then the following rules hold:
\begin{align*}
  \laplace_\I (g\cdot h) &= g \cdot\laplace_\I h+h \cdot\laplace_\I g+\Bigl(\partialmatrix{U} g\Bigr)^\trans\cdot\I\cdot\Bigl(\partialmatrix{U} h\Bigr)+\Bigl(\partialmatrix{U} h\Bigr)^\trans\cdot\I\cdot\Bigl(\partialmatrix{U} g\Bigr)\\
  \eulerop(g\cdot h) &= g\cdot\eulerop h+h\cdot\eulerop g
\end{align*}
\end{lemma}
We omit the proof as the claim follows by a straightforward calculation.
The part of (\ref{align_solutions}) that depends on the subspace of $\R^{m\times n}$, on which the quadratic form is negative definite, satisfies a slightly different system of partial differential equations than the one given in Lemma \ref{differentialequation_fr}, as an additional exponential factor occurs.
\begin{lemma}\label{satz_basis_negativ}
Let $\basis{\beta}{n}$ denote a basis of $\myhomopol{\beta}$.
We consider the system of partial differential equations
\begin{align}\label{align_differentialequation_fs}
\mathcal{D}_{-I}f=-(\beta +m)\cdot I\cdot f.
\end{align}
A finite basis of all solutions of \eqref{align_differentialequation_fs} that also satisfy the growth condition $f(U)\ko \exp(\pi \tr (U^\trans U))\in \mathcal{S}(\R^{m \times n})$ is given by
the functions
\begin{align*}
f_P(U):=\exp \bigl(-\tr \laplace/8\pi\bigr)\bigl(P(U)\bigr)\ko \exp\bigl(-2\pi \tr(U^\trans U)\bigr)\quad\text{with}\quad P\in \basis{\beta}{n}.
\end{align*}
\end{lemma}
\begin{proof}
We define $g(U):=\exp (-2\pi \tr(U^\trans U))$ and $h_P(U):=\exp (-\tr \laplace/8\pi)(P(U))$.
Both functions satisfy systems of partial differential equations similar to (\ref{align_differentialequation_fs}):
we check that
$$\bigl(\eulerop g\bigr)(U)=-4\pi g(U)\cdot U^\trans U$$ and
$$ \bigl(\laplace g\bigr)(U)=-4\pi m \cdot I\cdot g(U)+16\pi^2 g(U) \cdot U^\trans U$$
hold. Hence, we have
\begin{align}\label{dglexponentialfunktion}
\frac{1}{4\pi}\laplace g=-\eulerop g-m\cdot I\cdot g.
\end{align}
Due to Proposition \ref{lemma_equivalence} for $A=I$, the identity
\begin{align}\label{align_hp}
\frac{1}{4\pi}\laplace h_P=\eulerop h_P -\beta\cdot I\cdot h_P
\end{align}
holds if and only if $P\in \myhomopol{\beta}$. 
Using the multiplication rules from Lemma \ref{lemma_produktregel}, and applying \eqref{dglexponentialfunktion} and \eqref{align_hp} in the calculation of $\Delta f_{P}=\Delta(g \cdot h_{P})$, we obtain
\begin{align*}
\frac{1}{4\pi}\laplace f_P&=\frac{1}{4\pi}\biggl(g\cdot \laplace h_P+h_P\cdot \laplace g
+\Bigl(\partialmatrix{U} g\Bigr)^\trans \Bigl(\partialmatrix{U} h_P\Bigr)+\Bigl(\partialmatrix{U} h_P\Bigr)^\trans \Bigl(\partialmatrix{U} g\Bigr)\biggr)\\
&=g\cdot \bigl(\eulerop h_P -\beta\cdot I\cdot h_P \bigr)+h_P\cdot \bigl( -\eulerop g-m\cdot I\cdot g\bigr)-2 g\cdot \eulerop h_P\\
&=-(\beta+m)\cdot I\cdot f_P-\eulerop f_P,
\end{align*}
where we use in the second step that $\eulerop^\trans h_P=\eulerop h_P$ holds, since $h_P$ satisfies (\ref{align_hp}) and the Laplacian is symmetric.

Analogously, one can show that for any solution $f$ of the system \eqref{align_differentialequation_fs} of partial differential equations, the function $h(U)=f(U)\ko\exp \bigl( 2\pi \tr (U^\trans U)\bigr)$ satisfies $\mathcal{D}_I h = \beta\cdot I \cdot h$. Since 
\[h(U)\ko\exp \bigl( -\pi \tr (U^\trans U)\bigr)=
f(U)\ko \exp(\pi \tr(U^\trans U))\in \mathcal{S}(\R^{m\times n})\]
by assumption, we can apply Proposition \ref{lemma_equivalence}, which states that we can describe a finite basis for all functions $h$ by $h_P=\exp (-\tr \laplace/8\pi)(P(U))$ with $P\in \basis{\beta}{n}$.
Thus, the functions $f_P$ form a finite basis of the solutions of (\ref{align_differentialequation_fs}) that satisfy the aforementioned growth condition.
\end{proof}
In the next lemma, we show that the substitution of $U$ by $S^{-1}U$ leads to the desired system of partial differential equations that is associated with $A$.
\begin{lemma}\label{lemma_dgl_equivalence_normalized}
Let $S\in \R^{m\times m}$ such that $A=(S^{-1})^\trans \I S^{-1}$ and consider the functions $f,f[S^{-1}]:\R^{m\times n}\longrightarrow \R$, where $f[S^{-1}](U)=f(S^{-1}U)$. The function $f$ satisfies (\ref{align_dgl}) if and only if $f[S^{-1}]$ satisfies (\ref{align_dgl_A}).
\end{lemma}
\begin{proof}
Let $i,j\in \{1,\ldots ,n\}$. It suffices to calculate
\begin{align*}
(\laplace_A)_{ij}\bigl(f[S^{-1}](U)\bigr)&=\sum\limits_{a,b=1}^m\partialmatrix{U_{ai}}(A^{-1})_{ab}\partialmatrix{U_{bj}}\bigl(f(S^{-1}U)\bigr)\\
&=\sum\limits_{a,b,\mu,\nu=1}^m (S^{-1})_{\mu a}(A^{-1})_{ab}(S^{-1})_{\nu b}\frac{\partial^2f}{\partial U_{\mu i}\partial U_{\nu j}}(S^{-1}U)\\
&=\sum\limits_{\mu=1}^m\I_{\mu\mu}\frac{\partial^2f}{\partial U_{\mu i}\partial U_{\mu j}}(S^{-1}U)=\bigl((\laplace_\I)_{ij} f\bigr)(S^{-1}U)
\end{align*}
and
\begin{align*}
\eulerentry \bigl(f[S^{-1}](U)\bigr)&=\sum\limits_{d=1}^m U_{di}\partialmatrix{U_{dj}}\bigl(f(S^{-1}U)\bigr)=\sum\limits_{d,\nu=1}^m U_{di}(S^{-1})_{\nu d}\frac{\partial f}{\partial U_{\nu j}}(S^{-1}U)\\
&=\sum\limits_{\nu=1}^m (S^{-1}U)_{\nu i}\frac{\partial f}{\partial U_{\nu j}}(S^{-1}U)=(\eulerentry f)(S^{-1}U)
\end{align*}
to deduce the claim.
\end{proof}
\begin{satz}\label{corollary_indef}
Let $\basis{\alpha}{n}$ denote a basis of $\myhomopol{\alpha}$ and let $A\in \Z^{m\times m}$ denote a non-degenerate symmetric matrix  of signature $(r,s)$. As in Remark \ref{remark_decomp}, we write $A$ as the sum of a positive semi-definite matrix $\Apos$ and a negative semi-definite matrix $\Aneg$ and define $M:=\Apos -\Aneg$. 
The functions $$f(U)=f_{\alpha,\beta}(U)= \exp \bigl( -\laplacian_M/8\pi\bigr)\bigl(P(U)\bigr)\ko  \exp\bigl( 2\pi \tr (U^\trans \Aneg U)\bigr),$$
where $P\in \myhomopol{\alpha +\beta}$ is given as the product $P(U)=P_r(\Upos)\cdot P_s(\Uneg)$ with $P_r\in \basis{\alpha}{n}\subset\myhomopol{\alpha}$ and $P_s\in \basis{\beta}{n}\subset\myhomopol{\beta}$ for $(\alpha,\beta)\in \N_0^2$ such that $\alpha-\beta=\lambda+s$, form a (possibly infinite) basis for the space of solutions $f:\R^{m\times n}\longrightarrow\R$ of \eqref{align_dgl_A} that additionally satisfy the growth condition $$f(U)\ko \exp \bigl( -\pi \tr (U^\trans A U) \bigr)\in \mathcal{S}(\R^{m\times n}).$$.
\end{satz}
\begin{proof}
We consider the case $A=\I$. First we take $(\alpha,\beta)\in \N_0^2$ with $\alpha-\beta=\lambda+s$ to be fixed and show that $f=f_{\alpha,\beta}$ solves \eqref{align_dgl}. As the eigenvectors of $A$ form the canonical basis of $\R^m$, the polynomial $P$ splits as $P(U)=P_r(\ur)\cdot P_s(\us)$, where $\ur\in \R^{r\times n}$ consists of the first $r$ rows of $U$ and $\us\in \R^{s\times n}$ of the last $s$ rows of $U$. The exponential part of $f$ has the form $\exp\bigl(-2\pi \tr(\us^\trans \us)\bigr),$ so we can write $f:=f_r\cdot f_s$, where $f_r$ denotes the part dependent on $\ur$ and $f_s$ the part dependent on $\us$.
By Lemma \ref{lemma_produktregel}, we have $$\eulerop f=f_r\cdot\eulerop f_s+f_s\cdot\eulerop f_r.$$
The expression
\begin{align*}
\laplace_\I f=\laplace_\I (f_r\cdot f_s)=f_r\cdot\laplace_\I f_s+f_s\cdot\laplace_\I f_r
+\Bigl(\partialmatrix{U} f_r\Bigr)^\trans\cdot\I\cdot\Bigl(\partialmatrix{U} f_s\Bigr)+\Bigl(\partialmatrix{U} f_r\Bigr)^\trans\cdot\I\cdot\Bigl(\partialmatrix{U} f_s\Bigr)
\end{align*}
simplifies to $$\laplace_\I f=f_r\cdot\laplace_\I f_s+f_s\cdot\laplace_\I f_r,$$
since
\begin{align*}
\partialmatrix{U} f_r=\begin{pmatrix}
\partialmatrix{U_r}f_r\\ \mathrm{O}
\end{pmatrix}\quad\text{and}\quad \partialmatrix{U} f_s=\begin{pmatrix}
\mathrm{O}\\ \partialmatrix{U_s}f_s
\end{pmatrix}.
\end{align*}
These relations also show that we can write $\laplace_\I f_r=\laplace_{I_r}f_r$ and $\laplace_\I f_s=-\laplace_{I_s}f_s$. Then we consider the system of partial differential equations depending on the first $r$ rows of $U$, where $f_r$ corresponds to the function $f$ from Lemma \ref{differentialequation_fr}. Independent from that, consider the part depending on the last $s$ rows of $U$ and apply Lemma \ref{satz_basis_negativ} for $f_s$. Putting these results together, we obtain
\begin{align*}
\laplace_\I f&=4\pi \cdot\Bigl( f_r\cdot\bigl( \eulerop f_s+(\beta+s)\cdot I\cdot f_s\bigr)+f_s\cdot\bigl( \eulerop f_r-\alpha\cdot I\cdot f_r\bigr)\Bigr)\\
&=4\pi\cdot \bigl( \eulerop(f_r\cdot f_s)+(-\alpha+\beta +s)\cdot I\cdot (f_r\cdot f_s)\bigr)\\
&=4\pi\cdot\bigl( \eulerop f - (\alpha -\beta -s)\cdot I \cdot f\bigr),
\end{align*}
where $\alpha -\beta -s=\lambda$.

To show that these functions form a basis of all solutions, we employ a similar argument as in the proof of Proposition \ref{lemma_equivalence}. Again, we use Vignéras' result to show that the solutions $f$ of \eqref{align_dgl} have a certain form: We define the function $\widetilde{f}(U):=f(U)\ko \exp \bigl( -\pi \tr (U^\trans \I U) \bigr)$, which is a Schwartz function by assumption.
Furthermore, identify $\R^{m\times n}$ with $\R^{m n}$ by writing $U\in \R^{m\times n}$ as a column vector in $\R^{mn}$. As we have 
$$\tr (U^\trans \I U)=\sum_{\nu=1}^n\Bigl(\sum_{\mu=1}^{r} U_{\mu\nu}^2-\sum_{\mu=r+1}^{r+s} U_{\mu\nu}^2\Bigr),$$
which equals the normalized quadratic form of signature $(rn,sn)$ on $\R^{mn}$, we write $\widetilde{f}$ as $$\widetilde{f}(U)=f(U)\ko \exp\biggl(-\pi\sum_{\nu=1}^n\Bigl(\sum_{\mu=1}^{r} U_{\mu\nu}^2-\sum_{\mu=r+1}^{r+s} U_{\mu\nu}^2\Bigr)\biggr).$$ As an $\mathcal{L}^2(\R^{m n})$-function, $\widetilde{f}$ is given in an orthogonal basis of Hermite functions $H_{\boldsymbol{k}}$ in $m n$ variables in the form of $\widetilde{f}=\sum_{\boldsymbol{k}\in \N_0^{mn}} c_{\boldsymbol{k}}\ko  H_{\boldsymbol{k}}$ with $c_{\boldsymbol{k}}\in\R$.
Since $f$ is a solution of
\begin{align}\label{align_differentialequation_trace}
\tr(\mathcal{D}_{\I})f=\lambda n\cdot f\quad(\lambda\in\Z),
\end{align}
we restrict the possible basis elements that appear in the expansion of $\widetilde{f}$: $$\widetilde{f}=\sum\limits_{\substack{\boldsymbol{k}\in \N_0^{mn}\\ \epsilon (\boldsymbol{k})=n(\lambda+s)}} c_{\boldsymbol{k}}\ko  H_{\boldsymbol{k}},\quad\text{where}\quad\epsilon(\boldsymbol{k}):=\sum_{\nu=1}^n\Bigl(\sum_{\mu=1}^{r} k_{\mu\nu}-\sum_{\mu=r+1}^{r+s} k_{\mu\nu}\Bigr)$$
Thus, as a consequence of Vignéras' result for genus 1, any solution of (\ref{align_differentialequation_trace}) is given as a (possibly infinite) linear combination of functions
\begin{align*}
f_{\boldsymbol{k}}(U):=H_{\boldsymbol{k}}(U)\ko \exp \biggl(\pi\sum_{\nu=1}^n\Bigl(\sum_{\mu=1}^{r} U_{\mu\nu}^2-\sum_{\mu=r+1}^{r+s} U_{\mu\nu}^2\Bigr)\biggr),
\end{align*}
where the Hermite functions on $\R^{m n}$ (respectively $\R^{m\times n}$) are given as product of one-dimensional Hermite functions:
\begin{align*}
H_{\boldsymbol{k}}(U)=\prod_{\mu=1}^{m}\prod_{\nu=1}^{n}H_{k_{\mu\nu}}(U_{\mu\nu})=p(U_r)\ko q(U_s)\ko \exp \biggl( -\pi \sum_{\mu=1}^{m}\sum_{\nu=1}^n U_{\mu\nu}^2\biggr)
\end{align*}
with polynomials $p,q$, which are defined on $U_r, U_s$ respectively. Rewriting $f_{\boldsymbol{k}}$ as
\begin{align}\label{align_fk}
f_{\boldsymbol{k}}(U)&=p(U_r)\ko q(U_s)\ko \exp \Bigl( -2\pi \sum_{\mu=r+1}^{r+s}\sum_{\nu=1}^n U_{\mu\nu}^2\Bigr)
=p(U_r)\ko q(U_s)\ko \exp\bigl( -2\pi \tr (\us^\trans \us)\bigr),
\end{align}
each solution of \eqref{align_dgl} is given as a linear combination of functions of the form (\ref{align_fk}). The system of partial differential equations \eqref{align_dgl} is separable, i.\,e. can be broken into the part that depends on $\ur$ and the part that depends on $\us$. Likewise, $f_{\boldsymbol{k}}$ is given by a polynomial factor $p$ depending on $U_r$ and a factor of the form $q(U_s)\ko \exp\bigl( -2\pi \tr (\us^\trans \us)\bigr)$, where $q$ also denotes a polynomial. 
We can write $\mathcal{D}_{\I}=\mathcal{D}^r+\mathcal{D}^s$ such that the differential operator $\mathcal{D}^r$ vanishes if we apply it to a function on $\us$, and has the form $\mathcal{D}_{I_r}$ when applying it to a function that is defined on $\ur$. Analogously, $\mathcal{D}^s$ only depends on $\us$ and is of the form $\mathcal{D}_{-I_s}$ when applying it to functions on $\us$.
So we have $\mathcal{D}_{\I} f_{\boldsymbol{k}}=\lambda\cdot I\cdot f_{\boldsymbol{k}}$ with 
\begin{align*}
(\mathcal{D}_{\I} f_{\boldsymbol{k}})(U)
=q(U_s)\ko\exp\bigl( -2\pi \tr (\us^\trans \us)\bigr)\ko\mathcal{D}^r\bigl(p(\ur)\bigr)
+p(\ur)\ko\mathcal{D}^s\bigl(q(U_s)\ko\exp\bigl( -2\pi \tr (\us^\trans \us)\bigr)\bigr). 
\end{align*}
For $f_{\boldsymbol{k}}(U)\neq 0$ we divide by $f_{\boldsymbol{k}}$ and obtain for each entry of the system  of partial differential equations a sum of two partial differential equations that depend on different variables and therefore have to admit constant solutions. It follows that a function $f_{\boldsymbol{k}}$ solving \eqref{align_dgl} is given as the product described in (\ref{align_fk}) with the additional restriction that
\begin{multline*}
\frac{\mathcal{D}^r\bigl(p(\ur)\bigr)}{p(\ur)}=C_r\quad\text{and}\quad \frac{\mathcal{D}^s\bigl(q(U_s)\ko \exp\bigl( -2\pi \tr (\us^\trans \us)\bigr)\bigr)}{q(U_s)\ko \exp\bigl( -2\pi \tr (\us^\trans \us)\bigr)}=C_s\\ \text{with }C_r,C_s\in\R^{n\times n}\text{ and } C_r+C_s=\lambda \cdot I.
\end{multline*}
We show that $C_r=\alpha\cdot I$ for some $\alpha\in \N_0$ holds and thus $C_s=(\lambda-\alpha)\cdot I$. By applying the operator $\exp\bigl(\laplacian/8\pi\bigr)$ to $p(U_r)$, we can deduce analogously to the proof of Lemma \ref{differentialequation_fr} that $p(U_r)$ satisfies $$\mathcal{D}^r\bigl(p(\ur)\bigr)=C_r\cdot p(\ur)$$ if and only if the polynomial $P_r(\ur):=\exp\bigl(\laplacian/8\pi\bigr)\bigl(p(U_r)\bigr)$ satisfies $\eulerop P_r=C_r\cdot P_r$.
We have shown in Lemma \ref{lemma_konstante} that this system of partial differential equations admits polynomial solutions only if $C_r=\alpha \cdot I$ with $\alpha\in\N_0$.

Thus, every solution $f$ of \eqref{align_dgl} is described by basis elements $f_{\boldsymbol{k}}$ that consist of two factors that depend on different variables: $p$ solves the system of partial differential equations in Proposition \ref{lemma_equivalence}, where we have shown that these solutions can be described by a basis of homogeneous polynomials of degree $\alpha$.
Similarly, the function $q(U_s)\ko \exp\bigl( -2\pi \tr (\us^\trans \us)\bigr)$ solves the system of equations in Lemma \ref{satz_basis_negativ}, where we also described a basis of solutions using homogeneous polynomials of degree $\beta$.
We conclude that all solutions of \eqref{align_dgl} are described by the functions $f_{\alpha,\beta}$ defined above with $\alpha,\beta\in \N_0$ such that $\alpha-\beta=\lambda+s$. Thus, the basis consists of infinitely many elements if $r,s>0$ (i.\,e. when $A$ is indefinite) and finitely many otherwise (i.\,e. when $A$ is positive or negative definite).

We substitute $U\mapsto S^{-1}U$ and apply Lemma \ref{lemma_dgl_equivalence_normalized} to obtain the result for the system of partial differential equations \eqref{align_dgl_A}. Note that for every basis element $f_{\alpha,\beta}$ the polynomial $P$ splits as $P(U)=P_r(\Upos)\cdot P_s(\Uneg)$ by assumption.
\end{proof}
\section{Construction of theta series with modular transformation behavior}\label{sectiontransformation}
In this section, we construct Siegel theta series, which transform like modular forms of weight $m/2+\lambda$, arising from the functions that we considered in the last section as solutions of $\mathcal{D}_A f=\lambda\cdot I \cdot f$.
We explicitly determine the transformation behavior of the theta series with respect to $Z\mapsto Z+S$ (for a symmetric matrix $S\in \mymatrix{n}{\Z}$) and $Z\mapsto -Z^{-1}$.
To state the next lemma, in which we describe the transformation behavior of $\mythetac{H}{K}$ with respect to the first-mentioned transformation, we introduce the following notation for matrices:
\begin{defi}
(a) For $M\in \mymatrix{\mu}{\Z}$ we define $M_0\in \mymatrix{\mu}{\Z}$ by $(M_0)_{ij}=M_{ii}$ for $i=j$ and zero otherwise.\\
(b) We write $\mathrm{1}_{\mu \nu}$ for a matrix with $\mu$ rows and $\nu$ columns, whose entries are all equal to 1.
\end{defi}
\begin{lemma}\label{lemma_translation}
Let $S\in \mymatrix{n}{\Z}$ denote a symmetric matrix.
With respect to $Z\mapsto Z+S$, the theta series from Definition \ref{align_definition_theta} transforms as follows:
\begin{align*}
\mythetac{H}{K}(Z+S)=\e\bigl(- \tr (H^\trans AHS)/2- \tr (S_0\mathrm{1}_{nm}A_0H)/2\bigr)\ko \mythetac{H}{\widetilde{K}}(Z)
\end{align*}
with
\begin{align*}
\widetilde{K}:=K+HS+\frac12 A^{-1}A_0\mathrm{1}_{mn}S_0
\end{align*}
\end{lemma}
\begin{proof}
Write $U=H+R$ with $R\in\Z^{m\times n}$ such that
\begin{align*}
\tr (U^\trans AUS)=\tr \bigl(H^\trans AHS\bigr)+2\tr\bigl((HS)^\trans AR\bigr)+\tr\bigl(R^\trans ARS\bigr).
\end{align*}
It is straightforward to see that
\begin{align*}
\tr \bigl(H^\trans AHS\bigr)+2\tr\bigl((HS)^\trans AR\bigr)=-\tr\bigl(H^\trans AHS\bigr)+2\tr\bigl((HS)^\trans AU\bigr).
\end{align*}
As $A$ and $S$ both denote symmetric matrices and $x^2\equiv x\smod{2}$ for any $x\in\Z$, we have
\begin{align*}
\tr\bigl(R^\trans ARS\bigr)\equiv \sum\limits_{\nu=1}^n\sum\limits_{\mu=1}^m R_{\mu\nu}A_{\mu\mu}S_{\nu\nu}\smod{2}.
\end{align*}
To rewrite the expression on the right-hand side in terms of matrices, we introduce the matrix $\mathrm{1}_{nm}\in\Z^{n\times m}$ that only contains 1's as entries and obtain
\begin{align*}
\e\bigl(\tr (R^\trans ARS)/2\bigr)&=\e\bigl(\tr (S_0\mathrm{1}_{nm}A_0R)/2\bigr)\\
&=\e\Bigl(\tr \bigl(( A^{-1}A_0\mathrm{1}_{mn}S_0)^\trans AU\bigr)/2- \tr\bigl(S_0\mathrm{1}_{nm}A_0H\bigr)/2\Bigr).\qedhere
\end{align*}
\end{proof}
In the following two sections, we determine how the theta series  behaves under $Z\mapsto -Z^{-1}$. To put it briefly, we calculate the Fourier transform of the summand and then apply the Poisson summation formula. We define the Fourier transform associated with the matrix $A$:
\begin{defi}
Let $f:\R^{m\times n}\longrightarrow \C$ such that $f\in \mathcal{S}(\R^{m\times n})$. Then $\widehat{f}\in \mathcal{S}(\R^{m\times n})$ denotes the Fourier transform
\begin{align*}
\widehat{f}(V):=\myintegral f(U)\ko\e\bigl(\tr ( V^\trans A  U)\bigr)\ko dU
\end{align*}
with $dU$ the Euclidean volume element.
\end{defi}
Note that we do not take the standard definition of the Fourier transform as a unitary operator here, but rather we obtain the additional normalizing factor $|\det A|^{-n/2}$. Consequently, the Poisson summation formula has the form
\begin{align}\label{align_poisson}
\sum\limits_{U\in\Z^{m\times n}} f(U)=\sum\limits_{V\in A^{-1}\Z^{m\times n}} \widehat{f}(V).
\end{align}
In Section \ref{section_transformation_positive}, we consider theta series associated with positive definite quadratic forms and give a set of examples of non-holomorphic Siegel modular forms. We obtain those results by generalizing the set-up of Freitag \cite{Fre83}. In Section \ref{sectionindefinite}, we see that a similar construction also yields theta series associated with indefinite quadratic forms that transform like Siegel modular forms.
\subsection{Theta series for positive definite quadratic forms}\label{section_transformation_positive}
In this section, $p:\R^{m\times n}\longrightarrow \R$ is a polynomial and $A\in \mymatrix{m}{\Z}$ is a symmetric positive definite matrix.
Following Freitag \cite{Fre83}, we first examine the series 
\begin{align}\label{align_theta_freitag}
\sum_{U\in \Z^{m\times n}} p(U)\ko\e\bigl(\tr(U^\trans A UZ)/2\bigr)\quad (Z\in \H_n).
\end{align}
We consider the operator $$\tr \laplace_A =\sum\limits_{\nu=1}^n \sum\limits_{\mu,\rho =1}^m \frac{\partial}{\partial U_{\mu\nu}}\ko (A^{-1})_{\mu \rho}\ko \frac{\partial}{\partial U_{\rho\nu}}$$
and define
\begin{align*}
\exp (c\laplacian_A ) \bigl(p(U)\bigr):=\sum\limits_{k=0}^\infty \frac{c^k}{k !}(\laplacian_A)^k \bigl(p(U)\bigr)\quad (c\in \C).
\end{align*}
Since we are assuming that $p$ is a polynomial, this sum is finite.
\begin{lemma}
The following rules hold for $a,b,c\in \C$ and $M\in \mymatrix{m}{\C},\ko N\in \mymatrix{n}{\C}$:
\begin{align}
\exp(a\laplacian_A)\bigl(\exp(b\laplacian_A)(p(U))\bigr)&=\exp\bigl((a+b)\laplacian_A\bigr)\bigl(p(U)\bigr)\label{propertylaplace1}\\[2pt]
\exp(c\laplacian_A)\bigl(p(aU)\bigr)&=\bigl(\exp(a^2c\laplacian_A)p\bigr)(aU)\label{propertylaplace2}\\[2pt]
\exp(c\laplacian_A)\bigl(p(UN)\bigr)&=\bigl(\exp\bigl(c\tr(N\ko \laplace_A\ko  N^\trans)\bigr)p\bigr)(UN)\label{propertylaplace3}\\[2pt]
\exp(c\laplacian_A)\bigl(p(MU)\bigr)&=\biggl(\exp\Bigl(c\tr\Bigl( \bigl(\partialmatrix{U}\bigr)^\trans\ko MA^{-1}M^\trans\ko \partialmatrix{U}\Bigr)\Bigr) p\biggr)(MU)\label{propertylaplace4}
\end{align}
\end{lemma}
\begin{proof}
We derive Property (\ref{propertylaplace1}) by considering the Cauchy product for the absolutely convergent series $\sum_{k=0}^\infty \frac{1}{k !}(a\laplacian_A)^k$ and $\sum_{k=0}^\infty \frac{1}{k !}(b\laplacian_A)^k.$ The identity (\ref{propertylaplace2}) follows immediately from (\ref{propertylaplace3}), when we set $N:=a\cdot I \in \mymatrix{n}{\C}$. To show (\ref{propertylaplace3}) we consider $p(UN)$ and apply the Laplacian. We have
\begin{align*}
\partialmatrix{U_{\mu \nu}}\bigl(p(UN)\bigr)=\sum\limits_{\ell=1}^m \sum\limits_{i=1}^n \frac{\partial p}{\partial U_{\ell i}} (UN)\ko \frac{\partial (UN)_{\ell i }}{\partial U_{\mu \nu}}=\sum\limits_{i=1}^n N_{\nu i} \frac{\partial p}{\partial U_{\mu i}} (UN),
\end{align*}
since $$\frac{\partial (UN)_{\ell i}}{\partial U_{\mu \nu}}=\begin{cases}
N_{\nu i}&\text{ if }\ell=\mu,\\
0&\text{ otherwise}.
\end{cases}$$
With the same argument we obtain
\begin{align*}
\frac{\partial^2}{\partial U_{\mu \nu}\partial U_{\rho\nu}}\bigl(p(UN)\bigr)=\sum\limits_{d,e=1}^n N_{\nu d}N_{\nu e}\frac{\partial^2p}{\partial U_{\rho d}\partial U_{\mu e}}(UN)
\end{align*}
and therefore
\begin{align*}
\laplacian_A \bigl(p(UN)\bigr)&=\sum\limits_{\nu=1}^n \sum\limits_{\mu,\rho =1}^m \frac{\partial}{\partial U_{\mu\nu}}\ko (A^{-1})_{\mu \rho}\ko \frac{\partial}{\partial U_{\rho\nu}}\bigl(p(UN)\bigr)\\
&=\sum\limits_{\mu ,\rho =1}^m\sum\limits_{\nu,d,e=1}^n N_{\nu d}N_{\nu e} (A^{-1})_{\mu \rho} \frac{\partial^2p} {\partial U_{\rho d}\partial U_{\mu e}}(UN)\\
&=\biggl(\sum\limits_{\nu=1}^n\sum\limits_{\mu ,\rho=1}^m \Bigl(N\bigl(\partialmatrix{U}\bigr)^\trans\Bigr)_{\nu \mu} (A^{-1})_{\mu \rho}\Bigl(\partialmatrix{U}N^\trans\Bigr)_{\rho\nu}\ko  p\biggr)(UN)\\
&=\bigl(\tr (N\ko  \laplace_A\ko  N^\trans)p\bigr)(UN).
\end{align*}
Rewriting the Laplacian in the sum then gives (\ref{propertylaplace3}).
Analogously we obtain (\ref{propertylaplace4}).
\end{proof}
We calculate the Fourier transform of the summands in the series (\ref{align_theta_freitag}). To shorten the calculation, we apply the following result by Freitag \cite[p.\,158f.]{Fre83}, who considers Gauss transforms: we have
\begin{align}\label{align_freitag}
\myintegral p(U+V)\ko\exp\bigl(-\pi\tr(U^\trans U)\bigr)\ko dU=\exp\bigl( \laplacian/4\pi\bigr) \bigl(p(V)\bigr).
\end{align}
Note that Freitag uses the normalized Laplace operator $\laplacian=\tr\laplace_I$.
In the next lemma, we see that for arbitrary polynomials $p$, the functions $p(U)\ko \e\bigl( \tr ( U^\trans A UZ)/2\bigr)$ are not necessarily eigenfunctions with regard to the Fourier transform:
\begin{lemma}\label{lemma_fouriertransform_general_p}
Let $f:\R^{m\times n}\longrightarrow \C,\ko  f(U):=p(U)\ko \e\bigl( \tr ( U^\trans A UZ)/2\bigr).$
The Fourier transform of $f$ is
\begin{align*}
\widehat{f}(V)=\det A^{-n/2}\det (-iZ)^{-m/2}\e\bigl(- \tr (V^\trans A VZ^{-1})/2\bigr)\ko\Bigl(\exp\bigl(i\ko\tr( \laplace_A Z^{-1})/4\pi\bigr)p\Bigr)(-VZ^{-1}).
\end{align*}
\end{lemma}
\begin{proof}
We rewrite Freitag's result (\ref{align_freitag}) to obtain a form that is suitable for the calculation of the Fourier transform: We substitute $U$ by $U+iV$, and as we examine a holomorphic integrand in several complex variables, we apply the global residue theorem 
(i.\,e. instead of integrating over $U$ one can integrate over $U+iV$ without changing the integral) and obtain
\begin{align}\label{align_gauss_transform}
\begin{split}
\myintegral p(U)\ko &\exp \bigl(-\pi \tr (U^\trans U)+2\pi i \tr (V^\trans U)\bigr)\ko dU\\
&=\exp \bigl( -\pi \tr (V^\trans V)\bigr)\myintegral p(U+iV)\ko \exp \bigl(-\pi \tr (U^\trans U)\bigr)\ko dU\\ 
&=\exp \bigl( -\pi \tr (V^\trans V)\bigr)\ko \bigl(\exp(\laplacian/4\pi)p\bigr) (iV).
\end{split}
\end{align}
To determine
\begin{align*}
\widehat{f}(V)=\myintegral p(U)\ko \e\bigl(\tr (U^\trans AUZ)/2+ \tr (V^\trans A U)\bigr)\ko dU,
\end{align*}
we set $Z=iY$ and substitute $U$ by $A^{-1/2}UY^{-1/2}$ ($A$ and $Y$ are positive definite symmetric matrices, so the same holds for the inverses and uniquely determined square roots):
\begin{align*}
\myintegral p(U)\ko &\exp\bigl(-\pi\ko \tr (U^\trans AUY)+2\pi i\ko  \tr (V^\trans AU)\bigr)\ko dU\\
&=\begin{multlined}[t]\det A^{-n/2}\ko \det Y^{-m/2} \myintegral  p(A^{-1/2}UY^{-1/2}) \\
\cdot\exp\Bigl(-\pi\ko \tr (U^\trans U)+2\pi i\ko  \tr \bigl(( A^{1/2}VY^{-1/2})^\trans U\bigr)\Bigr)\ko dU\end{multlined}
\end{align*}
This is (\ref{align_gauss_transform}) evaluated at $A^{1/2}VY^{-1/2}$ with a slightly changed argument in the polynomial $p$. We apply (\ref{propertylaplace3}) and (\ref{propertylaplace4}) and use that $Y$ is symmetric to write $\widehat{f}$ as
\begin{align*}
\det A^{-n/2}\ko \det Y^{-m/2}  \exp\bigl(-\pi \tr (V^\trans A VY^{-1})\bigr)
\ko \Bigl(\exp \bigl(\tr(\laplace_A\ko  Y^{-1})/4\pi\bigr)p\Bigr)(iVY^{-1}).
\end{align*}
As the integrand is a holomorphic function, we resubstitute $Y=-iZ$ (for the inverse we have $Y^{-1}=iZ^{-1}$) and deduce the claim by analytic continuation.
\end{proof}
In order to obtain an eigenfunction under the Fourier transformation, Freitag \cite{Fre83} chooses $p$ to be a harmonic polynomial, i.\,e. $(\laplacian) p=0$ and $p(UN)=\det N^{\alpha} p(U)$ holds for all $N\in \mymatrix{n}{\C}$. We consider the more general class of polynomials
\begin{align}\label{align_choice_of_p}
p_Z(U):=\exp\bigl(-\tr(\laplace_A\ko Y^{-1})/8\pi\bigr)\bigl(P(U)\bigr)\quad \text{with}\quad P\in \myhomopol{\alpha}.
\end{align}
We described the vector space $P\in \myhomopol{\alpha}$ in Section \ref{section_homogeneity}, where we have also seen that the functions $p(U):=\exp(-\tr(\laplace_A)/8\pi)\bigl(P(U)\bigr)$ with $P\in \myhomopol{\alpha}$ form a basis for the vector space of the solutions of $\mathcal{D}_A f=\alpha\cdot I \cdot f$. 
The slightly modified functions in (\ref{align_choice_of_p}) depend on the imaginary part $Y$ of $Z$, which means that we lose holomorphicity in the construction of the theta series.
However, for harmonic polynomials $P$, we obtain the holomorphic theta series considered by Freitag. Note that this is basically a generalization of Borcherds' construction for $n=1$ in \cite{Bor98}, see Remark \ref{remark_borcherds} for a more detailed explanation.
\begin{lemma}\label{lemma_ft_hom_pos}
Let $p_Z$ denote a polynomial from (\ref{align_choice_of_p}) and define
$$f_Z(U):=p_Z(U)\ko  \e\bigl( \tr ( U^\trans A UZ)/2\bigr).$$ The Fourier transform of $f_Z$ is
\begin{align*}
\widehat{f_Z}(V)=i^{-mn/2}\det A^{-n/2} \det(-Z^{-1})^ {m/2+\alpha}\ko f_{-Z^{-1}}(V).
\end{align*}
\end{lemma}
\begin{proof}
We apply Lemma \ref{lemma_fouriertransform_general_p} and then use the linearity of the trace and Property (\ref{propertylaplace1}):
\begin{align}\label{align_fouriertransform_fz}
\begin{split}
\myintegral p_Z(U)\ko &\e\bigl(\tr (U^\trans AUZ)/2+ \tr (V^\trans A U)\bigr)\ko dU\\
&=\begin{multlined}[t] 
\det A^{-n/2}\ko \det (-iZ)^{-m/2} \e\bigl(- \tr (V^\trans A  VZ^{-1})/2\bigr)\\
\cdot\left(\exp\bigl( i\ko\tr\bigl(\laplace_A\ko  Z^{-1}\bigr)/4\pi-\tr( \laplace_A\ko  Y^{-1})/8\pi\bigr)P\right)(-VZ^{-1})\end{multlined}\\[2pt]
&=\begin{multlined}[t]
\det A^{-n/2}\ko \det (-iZ)^{-m/2} \e\bigl(- \tr (V^\trans A VZ^{-1})/2\bigr)\\
\cdot\left(\exp\bigl(-\tr ( \laplace_A\ko  (Y^{-1}-2iZ^{-1}))/8\pi\bigr)P\right)(-VZ^{-1})\end{multlined}
\end{split}
\end{align}
If $\widetilde{Y}$ denotes the imaginary part of $-Z^{-1}$, the identity $\widetilde{Y}=\overline{Z}^{-1} Y Z^{-1}$ holds by (\ref{imaginary_part}). Hence,
\begin{align*}
Y^{-1}-2iZ^{-1}=Y^{-1}(Z-2iY)Z^{-1}=Y^{-1}\overline{Z} Z^{-1}=Z^{-1}(ZY^{-1}\overline{Z})Z^{-1}=Z^{-1}\widetilde{Y}^{-1} Z^{-1}.
\end{align*}
The matrix $Z$ is symmetric and therefore also its inverse $Z^{-1}$, which means that we can rewrite (\ref{align_fouriertransform_fz}) as follows:
\begin{multline*}
\det A^{-n/2}\ko \det (-iZ)^{-m/2} \e\bigl(- \tr (V^\trans A VZ^{-1})/2\bigr)\\
\cdot\left(\exp\Bigl( -\tr \bigl((-Z^{-1})\laplace_A(-Z^{-1})^\trans\ko \widetilde{Y}^{-1}\ko  \bigr)/8\pi\Bigr)P\right)(-VZ^{-1})
\end{multline*}
Using Property (\ref{propertylaplace3}) and the homogeneity of $P\in \myhomopol{\alpha}$, we conclude that the Fourier transform of $f_Z$ has the form
\begin{align*}
\widehat{f_Z}(V)&=
\det A^{-n/2}\ko \det (-iZ)^{-m/2}\ko  \det(-Z^{-1})^{\alpha}\ko  \e\bigl(- \tr (V^\trans A  VZ^{-1})/2\bigr)
\ko\Bigl(\exp\bigl(-\tr (\laplace_A\ko  \widetilde{Y}^{-1})/8\pi\bigr)P\Bigr)(V) \\
&=\det A^{-n/2}\ko \det (-iZ)^{-m/2}\ko  \det(-Z^{-1})^{\alpha}\ko f_{-Z^{-1}}(V).
\end{align*}
Separating constant factors and factors depending on the determinants of $A$ and $Z$, we deduce the claim.
\end{proof}
This construction yields theta series that transform like Siegel modular forms: 
\begin{satz}\label{proposition_transformation_posdef}
Let $A\in \Z^{m\times m}$ denote a positive definite symmetric matrix and $p$ the polynomial defined as $p(U)=\exp\bigl(-\laplacian_A/8\pi\bigr)\bigl(P(U)\bigr)$ with $P\in \myhomopol{\alpha}$. For the corresponding theta series $\mythetac{H}{K}$ given in Definition \ref{align_definition_theta} we have
\begin{align*}
\mythetac{H}{K}(-Z^{-1})=i^{-mn/2}\det A^{-n/2}\det Z^{m/2+\alpha}\e \bigl( \tr (H^\trans A K)\bigr)
\ko\sum\limits_{J\in A^{-1}\Z^{m\times n} \operatorname{mod} \Z^{m\times n} } \mythetac{J+K}{-H}(Z).
\end{align*}
\end{satz}
\begin{proof}
We recall the definition of $\mythetac{H}{K}$, which is
\begin{align*}
\mythetac{H}{K}(Z)=\det Y^{-\alpha/2}\sum\limits_{U\in H+\Z^{m\times n}}p(UY^{1/2})\ko \e\bigl( \tr(U^\trans AUZ)/2+ \tr(K^\trans AU)\bigr).
\end{align*}
We use Property (\ref{propertylaplace3}) and the homogeneity property of $P$ to rewrite $$\det Y^{-\alpha /2}\ko p(UY^{1/2})=\exp \bigl(-\tr(\laplace_A Y^{-1})/8\pi\bigr)\bigl(P(U)\bigr)=p_Z(U),$$
and analogously $p_{-Z^{-1}}(U)=\det \widetilde{Y}^{-\alpha/2}\ko p(U\widetilde{Y}^{1/2})$.
That means the theta series has the form
\begin{multline*}
\mythetac{H}{K}(-Z^{-1})=\sum\limits_{U\in \Z^{m\times n}}\Big\lbrace p_{-Z^{-1}}(U+H)\\
\cdot \e\bigl(- \tr((U+H)^\trans A(U+H)Z^{-1})/2+ \tr(K^\trans A(U+H))\bigr)\Big\rbrace.
\end{multline*}
By Lemma \ref{lemma_ft_hom_pos}, the Fourier transform of the summand equals
\begin{align*}
i^{-mn/2}\det A^{-n/2}\det Z^{m/2+\alpha}\ko p_Z(V+K)
\ko\e \bigl(\tr ((V+K)^\trans A(V+K)Z)/2- \tr (H^\trans A V)\bigr).
\end{align*}
The summands in the theta series are Schwartz functions as $A$ denotes a positive definite quadratic form. Hence, we apply the Poisson summation formula \eqref{align_poisson}, and obtain
\begin{align*}
\mythetac{H}{K}(-Z^{-1})
&=\begin{multlined}[t]
i^{-mn/2}\det A^{-n/2}\det Z^{m/2+\alpha}\e\bigl( \tr (H^\trans AK)\bigr)\\
\cdot \sum\limits_{V\in K+A^{-1}\Z^{m\times n}}p_Z(V)\ko \e\bigl(\tr(V^\trans AVZ)/2- \tr(H^\trans AV)\bigr)\end{multlined}\\
&=\begin{multlined}[t]
i^{-mn/2}\det A^{-n/2}\det Z^{m/2+\alpha}\e\bigl( \tr (H^\trans AK)\bigr)\\
\cdot\sum\limits_{J\in A^{-1}\Z^{m\times n} \operatorname{mod} \Z^{m\times n} } \mythetac{J+K}{-H}(Z),\end{multlined}
\end{align*}
which completes the proof.
\end{proof}
\begin{bsp}\label{example_positive}
For $m\equiv 0\smod 8$, we choose an even unimodular matrix $A\in \Z^{m\times m}$, which means in particular that $\det A=1$ and $A^{-1}\in \Z^{m\times m}$. Considering the theta series $$\mythetac{\mathrm{O}}{\mathrm{O}}(Z)=\det Y^{-\alpha/2}\sum\limits_{U\in \Z^{m\times n}}p(UY^{1/2})\ko \e\bigl(\tr(U^\trans AUZ)/2\bigr),$$ we have $\mythetac{\mathrm{O}}{\mathrm{O}}(Z+S)=\mythetac{\mathrm{O}}{\mathrm{O}}(Z)$ for any symmetric matrix $S \in \mymatrix{n}{\Z}$ by Lemma \ref{lemma_translation} and $\mythetac{\mathrm{O}}{\mathrm{O}}(-Z^{-1})=\det Z^{m/2+\alpha}\ko \mythetac{\mathrm{O}}{\mathrm{O}}(Z)$ for a polynomial $p$ as chosen in Proposition \ref{proposition_transformation_posdef}.
Thus, $\mythetac{\mathrm{O}}{\mathrm{O}}$ is a non-holomorphic Siegel modular form of weight $m/2+\alpha$ on the full Siegel modular group $\Gamma_n$.
\end{bsp}
\subsection{Theta series for indefinite quadratic forms}\label{sectionindefinite}
In this section, we consider theta series associated with non-degenerate symmetric matrices $A\in \mymatrix{m}{\Z}$ with signature $(r,s)$, where $s\geq 0$. As described in Remark \ref{remark_decomp}, we decompose $A=\Apos +\Aneg$ by employing the matrix of normalized eigenvectors $S\in\R^{m\times m}$ so that we obtain the associated majorant matrix $\Amaj=\Apos -\Aneg$ and the projections $U^{\pm}$ of $U$ into the positive and negative subspaces of $\R^{m\times n}$ respectively.
We replace the polynomials $p$ that were defined in Proposition \ref{proposition_transformation_posdef} by functions of the form
\begin{align}\label{align_coefficient_indefinite}
g(U):=\exp\bigl( -\laplacian_{\Amaj}/8\pi\bigr)\bigl(P(U)\bigr)\ko \exp\bigl(2\pi \tr(U^\trans \Aneg U )\bigr),
\end{align}
where $P\in \myhomopol{\alpha +\beta}$ is given as the product $P(U)=P_{\alpha}(\Upos)\cdot P_{\beta}(\Uneg)$ with $P_{\alpha}\in \myhomopol{\alpha}$ and $P_{\beta}\in \myhomopol{\beta}$. For $\alpha-\beta=\lambda+s$, we know from Section \ref{section_coeff_diff_eq_pos} that these functions are solutions of $\mathcal{D}_A f=\lambda\cdot I\cdot f$. Of course, we can also replace $g$ by a linear combination of functions of this type, under the assumption that $(\alpha,\beta)\in \N_0^2$ such that $\alpha-\beta=\lambda+s$, to construct modular Siegel theta series. 
However, it is sufficient for the proof of Theorem \ref{maintheorem} and simplifies the following calculations just to consider $g$ as defined above, since these functions in particular include the basis elements of the vector space of solutions of $\mathcal{D}_A f=\lambda\cdot I\cdot f$.

In analogy with the last section, we define
\begin{align*}
g_Z(U):=\exp\bigl( -\tr (\laplace_{\Amaj} Y^{-1})/8\pi\bigr)\bigl(P(U)\bigr)\ko \exp\bigl(2\pi \tr(U^\trans \Aneg U Y )\bigr)
\end{align*}
and
\begin{align*}
f_Z(U):=g_Z(U)\ko  \e\bigl( \tr ( U^\trans A  UZ)/2\bigr).
\end{align*}
For $s=0$, we get back the functions from Lemma \ref{lemma_ft_hom_pos}, so we use the same notation.
\begin{lemma}\label{lemma_ft_hom_indef}
The Fourier transform of
$f_Z(U)=g_Z(U)\ko  \e\bigl( \tr ( U^\trans A UZ)/2\bigr)$
is
\begin{align*}
\widehat{f_Z}(V)=i^{-mn/2} (-1)^{\beta s} |\det A|^{-n/2} \det (-Z^{-1})^{r/2+\alpha}&\det \overline{Z}^{-(s/2+\beta)}\ko f_{-Z^{-1}}( V).
\end{align*}
\end{lemma}
\begin{proof}
We change the basis of $\R^m$ by the substitution of $U\mapsto SU$ to obtain a part that depends on the first $r$ rows of $U$ (again, we denote this part of the matrix by $\ur$) and one part that depends on the last $s$ rows of $U$ (analogously, we denote this part by $\us$):
\begin{align*}
\begin{split}
\myintegral g_Z(U)\ko  &\e\bigl( \tr ( U^\trans A UZ)/2+\tr (U^\trans AV)\bigr)\ko dU\\
&=\begin{multlined}[t]\myintegral \exp\bigl( -\tr (\laplace_{\Amaj} Y^{-1})/8\pi\bigr)\bigl(P(U)\bigr)\\
\cdot\exp\bigl(2\pi \tr(U^\trans \Aneg U Y )+\pi i\tr(U^\trans AUZ)+2\pi i\tr(U^\trans AV)\bigr)\ko dU\end{multlined}\\
&=\begin{multlined}[t]\det S^n\myintegral  \exp\bigl( -\tr (\laplace Y^{-1})/8\pi\bigr)\bigl(P(SU)\bigr)\\
\cdot\exp\bigl(-2\pi \tr(\us^\trans \us Y )+\pi i\tr(U^\trans\I UZ)+2\pi i\tr(U^\trans \I S^{-1} V)\bigr)\ko dU\end{multlined}
\end{split}
\end{align*}
We can split up the integral, as the polynomial $P$ factors as a polynomial dependent on $\ur$ and $\us$ respectively. We now apply the results for positive definite quadratic forms. By  Lemma \ref{lemma_ft_hom_pos}, we obtain
\begin{align}\label{align_fouriertransformation_ur}
\begin{split}
\int\limits_{\R^{r\times n}}\exp \biggl(-\frac{1}{8\pi}&\tr \Bigl(\Bigl( \partialmatrix{\ur}\Bigr)^\trans\partialmatrix{\ur}\ko Y^{-1}\Bigr)\biggr)\bigl(P_{\alpha}(\ur)\bigr)\ko\e\bigl(\tr(\ur^\trans \ur Z)/2+\tr(\vr^\trans \ur)\bigr)\ko d\ur\\[2pt]
&=\begin{multlined}[t] i^{-rn/2} \det (-Z^{-1})^{r/2+\alpha}\ko \e\bigl(-\tr  (\vr^\trans   \vr Z^{-1})/2\bigr)\\[2pt]
\cdot\exp\biggl(-\frac{1}{8\pi}\tr \Bigl(\Bigl(\partialmatrix{\vr}\Bigr)^\trans\partialmatrix{\vr}\ko \widetilde{Y}^{-1}\Bigr)\biggr)\bigl(P_{\alpha}(\vr)\bigr).\end{multlined}
\end{split}
\end{align}
We treat the part that depends on the negative definite subspace like an expression that is associated with a positive definite quadratic form given by $I_s$ and consider $-\overline{Z}\in \H_n$ as variable in the Siegel upper half-space. Also note that by (\ref{imaginary_part}) we have
\begin{align}\label{align_imaginary}
\overline{Z}^{-1}=\overline{Z}^{-1}ZZ^{-1}=\overline{Z}^{-1}(\overline{Z}+2iY)Z^{-1}=Z^{-1}+2i\overline{Z}^{-1}YZ^{-1}=Z^{-1}+2i\widetilde{Y}.
\end{align}
In particular, $\Im \bigl(\overline{Z}^{-1}\bigr)=\Im (-Z^{-1})=\widetilde{Y}$ and thus we have
\begin{align}\label{align_fouriertransformation_us}
\begin{split}
\int\limits_{\R^{s\times n}}\exp \biggl(-\frac{1}{8\pi}&\tr \Bigl(\Bigl(\partialmatrix{\us}\Bigr)^\trans\partialmatrix{\us}\ko Y^{-1}\Bigr)\biggr)\bigl(P_{\beta}(\us)\bigr)\ko\e\bigl(- \tr(\us^\trans \us\overline{Z})/2-\tr(\vs^\trans \us)\bigr)\ko d\us\\
&=\begin{multlined}[t]i^{-sn/2}(-1)^{\beta s}\det (\overline{Z}^{-1})^{s/2+\beta}\ko \e\bigl( \tr  (\vs^\trans   \vs \overline{Z}^{-1})/2\bigr)\\[2pt]
\cdot\exp\biggl(-\frac{1}{8\pi}\tr\Bigl(\Bigl(\partialmatrix{\vs}\Bigr)^\trans\partialmatrix{\vs}\widetilde{Y}^{-1}\Bigr)\biggr)\bigl(P_{\beta}(\vs)\bigr),\end{multlined}
\end{split}
\end{align}
where we evaluate the Fourier transform for $-\vs$, and use Property (\ref{propertylaplace3}) and the identity $P_{\beta}(-\vs)=(-1)^{\beta s} P_{\beta}(\vs)$ to rewrite the expression.
We now consider the product of (\ref{align_fouriertransformation_ur}) and (\ref{align_fouriertransformation_us}) (we make use of (\ref{align_imaginary}) again to rewrite the exponential factor) and obtain:
\begin{align*}
\myintegral \hspace{-2.4pt}&\exp\bigl( -\tr (\laplace Y^{-1})/8\pi\bigr)(P(SU))\ko
\exp\bigl(-2\pi \tr(\us^\trans \us Y )+\pi i\tr(U^\trans\I UZ)+2\pi i\tr(U^\trans\I V)\bigr)\ko dU\\
&=\begin{multlined}[t]i^{-mn/2} (-1)^{\beta s} \det (-Z^{-1})^{r/2+\alpha}\det \overline{Z}^{-(s/2+\beta)}\ko  \exp\bigl( -(\laplacian\widetilde{Y}^{-1})/8\pi\bigr)(P(SV))\\[2pt]
\cdot\exp\bigl(-\pi i \tr  (V^\trans  \I V Z^{-1})-2\pi  \tr  (\vs^\trans   \vs \widetilde{Y})\bigr)\end{multlined}
\end{align*}
We evaluate this integral at $S^{-1}V$ to complete the proof.
Without loss of generality, we can assume that $\det S>0$ and therefore write $\det S^n$ as $|\det A|^{n/2}$.
\end{proof}
Thus, we can state a more general version of Proposition \ref{proposition_transformation_posdef} for Siegel theta series for indefinite quadratic forms.
\begin{satz}\label{proposition_transformation_indef}
Let $\lambda=\alpha-\beta-s$ and let $g:\R^{m\times n}\longrightarrow \R$ define a function from (\ref{align_coefficient_indefinite}). 
The theta series of the form
\begin{align*}
\mythetac{H}{K}(Z)=\det Y^{-\lambda/2}\sum\limits_{U\in H+\Z^{m\times n}} g(UY^{1/2})\ko \e\bigl( \tr(U^\trans A UZ)/2+ \tr (K^\trans AU)\bigr)
\end{align*}
transforms as follows:
\begin{multline*}
\mythetac{H}{K}(-Z^{-1})=i^{-mn/2}(-1)^{(s/2+\beta)n+\beta s}|\det A|^{-n/2}\det Z^{(r-s)/2+\alpha-\beta}
\ko\e \bigl( \tr (H^\trans A K)\bigr)\\ \cdot\sum\limits_{J\in A^{-1}\Z^{m\times n} \operatorname{mod} \Z^{m\times n} } \mythetac{J+K}{-H}(Z)
\end{multline*}
\end{satz}
\begin{proof}
We use the same approach as in the proof of Proposition \ref{proposition_transformation_posdef}. By Property \eqref{propertylaplace3}, we have $g_{-Z^{-1}}(U)=\det \widetilde{Y}^{-(\alpha+\beta)/2}g(U\widetilde{Y}^{1/2})$, and thus we rewrite the theta series as
\begin{multline*}
\mythetac{H}{K}(-Z^{-1})=\det \widetilde{Y}^{s/2+\beta}\sum\limits_{U\in \Z^{m\times n}}\Big\lbrace g_{-Z^{-1}}(U+H)\\ 
\cdot\e\bigl(- \tr ( (U+H)^\trans A (U+H)Z^{-1})/2+ \tr(K^\trans A(U+H))\bigr)\Big\rbrace.
\end{multline*}
By Lemma \ref{lemma_ft_hom_indef}, the Fourier transform of the summand equals
\begin{multline*}
i^{-mn/2}(-1)^{\beta s}|\det A|^{-n/2}\det Z^{r/2+\alpha}\det (-\overline{Z})^{s/2+\beta}g_Z(V+K)\\
\cdot\e \bigl(\tr ((V+K)^\trans A(V+K)Z)/2- \tr (H^\trans A V)\bigr).
\end{multline*}
By \eqref{align_imaginary}, we have $\widetilde{Y}=\overline{Z}^{-1}YZ^{-1}$ and thus rewrite
\begin{align*}
\det \widetilde{Y}^{s/2+\beta}\det Z^{r/2+\alpha}\det (-\overline{Z})^{s/2+\beta}=
(-1)^{(s/2+\beta)n} \det Z^{(r-s)/2+\alpha-\beta}\det Y^{s/2+\beta}.
\end{align*}
As $\det Y^{s/2+\beta} g_Z(U)=\det Y^{-\lambda/2}g(UY^{1/2})$ by Property (\ref{propertylaplace3}), we have
\begin{multline*}
\mythetac{H}{K}(-Z^{-1})=i^{-mn/2}(-1)^{(s/2+\beta)n+\beta s} |\det A|^{-n/2} \det Z^{(r-s)/2+\alpha-\beta}\e\bigl(\tr(H^\trans AK)\bigr)\\
\cdot \det Y^{-\lambda/2}\sum\limits_{V\in K+A^{-1}\Z^{m\times n}} g(VY^{1/2})\e \bigl(\tr(V^\trans AVZ)/2-\tr (V^\trans AH)\bigr).
\end{multline*}
Again, we write
\begin{multline*}
\det Y^{-\lambda/2}\sum\limits_{V\in K+A^{-1}\Z^{m\times n}} g(VY^{1/2})\ko \e \bigl( \tr(V^\trans AVZ)/2-\tr (V^\trans AH)\bigr)\\
=\sum\limits_{J\in A^{-1}\Z^{m\times n} \operatorname{mod} \Z^{m\times n} } \mythetac{J+K}{-H}(Z),
\end{multline*}
which completes the proof.
\end{proof}
\begin{bsp}\label{example_negative}
We obtain examples of non-holomorphic Siegel modular forms on the full Siegel modular group if $H=K=\mathrm{O}$ and $A$ is an even unimodular matrix and additionally $i^{mn/2}\ko (-1)^{(s/2+\beta)n+\beta s}=1$ holds.
Note that an even symmetric unimodular matrix of indefinite signature $(r,s)$ only exists when $r-s\equiv 0\smod{8}$ and is isomorphic to $H_2^k\oplus (\pm E_8)^{\ell}$ with $k=\min\{r,s\}$ and $\ell=|r-s|/8$, where $H_2=\left(\begin{smallmatrix}
0&1\\1&0
\end{smallmatrix}\right)$ and $E_8$ represents the equivalence class of all even unimodular positive definite matrices of rank 8 (note that we take $E_8$ if $r>s$ and $-E_8$ if $r<s$), see for example Husemoller and Milnor \cite[p.\,24-26]{HM73} for more details.
\end{bsp}
\subsection{Proof of Theorem \ref{maintheorem}}
In Section \ref{section_vigneras}, we introduced the $n\times n$ system of partial differential equations $\mathcal{D}_A f=\lambda\cdot I\cdot f$ and determined a basis for all the solutions $f$ that additionally satisfy the growth condition $f(U)\ko \exp(-\pi \tr(U^\trans AU))\in \mathcal{S}(\R^{m\times n})$ (see Proposition \ref{corollary_indef}). 
In this section, we have determined the modular transformation behavior of the associated Siegel theta series $\theta_{H,K,f,A}$ by explicitly calculating the transformation formulas for the generators $Z\mapsto Z+S$  (see Lemma \ref{lemma_translation}) and $Z\mapsto -Z^{-1}$ (see Proposition \ref{proposition_transformation_indef}) of the Siegel modular group. For an even matrix $A$ and $\lambda=\alpha-\beta-s$, the theta series $\theta_{\mathrm{O},\mathrm{O},f,A}$ transforms like a Siegel modular form of genus $n$ and weight $m/2+\lambda$ on some congruence subgroup of $\Gamma_n$.
This proves Theorem \ref{maintheorem}.

\begin{acknowledgments}
The author is grateful to Sander Zwegers for suggesting the topic and lots of helpful advice. Besides, the author would like to kindly thank Jens Funke for useful remarks during the DMV Jahrestagung 2019, Markus Schwagenscheidt for pointing out how the system of differential equations can be reduced in order to apply Vignéras' result, Ben Wright for carefully checking the manuscript for linguistic deficiencies, and the reviewers who contributed many insightful comments that improved the exposition of this paper significantly.
\end{acknowledgments}
\bibliography{References}

\end{document}